\documentclass[11pt,twoside,a4paper]{article}
\title{} \author{} \date{}
\usepackage{latexsym,amssymb,times}
\input amssym.def
\newtheorem{te}{Theorem}[section]
\newtheorem{prop}[te]{Proposition}

\newtheorem{fac}[te]{Fact}
\newtheorem{cla}[te]{Claim}

\newtheorem{df}[te]{Definition}
\newtheorem{rem}[te]{Remark}

\def\dok{\noindent{\bf Proof. }}
\def\kdok{\hfill $\Box$ \par \vspace*{2mm} }
\def\a{\alpha}
\def\b{\beta}
\def\d{\delta}
\def\f{\varphi}

\def\o{\omega}
\def\k{\kappa}
\def\l{\lambda}
\def\r{\rho}

\def\t{\tau}
\def\ve{\varepsilon}
\def\S{{\mathbb S}}
\def\B{{\mathbb B}}
\def\N{{\mathbb N}}
\def\X{{\mathbb X}}
\def\Y{{\mathbb Y}}
\def\Z{{\mathbb Z}}
\def\A{{\mathbb A}}
\def\C{{\mathbb C}}

\def\BI{{\mathbb I}}
\def\BE{{\mathbb E}}
\def\CP{{\mathcal P}}
\def\CC{{\mathcal C}}
\def\CT{{\mathcal T}}
\def\c{{\mathfrak{c}}}
\def\la{\langle}
\def\ra{\rangle}
\def\id{\mathop{\mathrm{id}}\nolimits}

\def\Sym{\mathop{\rm Sym}\nolimits}
\def\Form{\mathop{\rm Form}\nolimits}
\def\Mod{\mathop{\rm Mod}\nolimits}
\def\Th{\mathop{\rm Th}\nolimits}
\def\bcd{\dot{\bigcup}}
\def\du{\mathrel{\dot{\cup}}}
\def\Perf{\mathop{\rm Perf}\nolimits}
\def\ar{\mathop{\rm ar}\nolimits}
\def\Clop{\mathop{\rm Clop}\nolimits}
\begin{document}
\thispagestyle{plain}
\begin{center}
           {\large \bf \uppercase{Vaught's conjecture for unions of products of rooted trees}}
\end{center}
\begin{center}
{\bf Milo\v s S.\ Kurili\'c}\footnote{Department of Mathematics and Informatics, Faculty of Sciences, University of Novi Sad,
                                      Trg Dositeja Obradovi\'ca 4, 21000 Novi Sad, Serbia.
                                      e-mail: milos@dmi.uns.ac.rs}
\end{center}
\begin{abstract}
\noindent
Let ${\mathcal C} ^{\rm rt}$ be the class of rooted trees and $\langle {\mathcal C} ^{\rm rt}\rangle _{\dot{\cup }\Pi}$ the minimal class of partial orders
containing ${\mathcal C} ^{\rm rt}$ and closed under isomorphism, finite direct products and finite disjoint unions.
Each poset from $\langle {\mathcal C} ^{\rm rt}\rangle _{\dot{\cup }\Pi}$ is isomorphic to one of the form ${\mathbb X}= \bcd _{i<n}\prod _{j<m_i}{\mathbb X}_i^j$,
where the factors ${\mathbb X}_i^j$ are rooted trees,
and defining  ${\mathcal T}=\mathop{\rm Th} ({\mathbb X})$,  ${\mathcal T} _i ^j=\mathop{\rm Th}({\mathbb X}_i^j)$, for $i<n$ and $j<m_i$, 
and $\kappa = \prod _{i<n}\prod _{j<m_i}I({\mathcal T} _i^j)$, we have
\begin{itemize}\itemsep -0.2mm
\item[\rm (a)] Vaught's conjecture is true for ${\mathcal T}$:  $I({\mathcal T})=\kappa $, if $\kappa\in \{ 1,\omega ,{\mathfrak{c}}\}$, and, otherwise, $I({\mathcal T}) \in [3,\omega)$;
\item[\rm (b)] ${\mathbb Y} \equiv {\mathbb X}$ iff $\;{\mathbb Y} \cong \dot{\bigcup}_{i<n}\prod _{j<m_i}{\mathbb Y} _i^j$, where ${\mathbb Y}_i^j\equiv {\mathbb X}_i^j$, for $i<n$ and $j<m_i$;
\item[\rm (c)] ${\mathbb E}\preccurlyeq {\mathbb X}$ iff $\;{\mathbb E} =\dot{\bigcup}_{i<n}\prod _{j<m_i}{\mathbb E}_i^j$, where ${\mathbb E}_i^j\preccurlyeq {\mathbb X}_i^j$, for $i<n$ and $j<m_i$;
\item[\rm (d)] ${\mathcal T}$ is atomic iff $\;{\mathcal T} _i^j$, for $i<n$ and $j<m_i$, are atomic;
               then $\dot{\bigcup}_{i<n}\prod _{j<m_i}{\mathbb A}_i^j$ is a countable atomic model of ${\mathcal T}$,
               where ${\mathbb A}_i^j$ is a countable atomic model of ${\mathcal T} _i^j$, for $i<n$ and $j<m_i$;
\item[\rm (e)] ${\mathcal T}$ is small iff $\;{\mathcal T} _i^j$, for $i<n$ and $j<m_i$, are small;
               then $\dot{\bigcup}_{i<n}\prod _{j<m_i}{\mathbb S}_i^j$ is a countably saturated model of ${\mathcal T}$,
               where ${\mathbb S}_i^j$ is a countably saturated model of ${\mathcal T}_i^j$, for $i<n$ and $j<m_i$.
\end{itemize}

\end{abstract}
\section{Introduction}\label{S1}
Vaught's conjecture, stated by Robert Vaught \cite{Vau},
says that the number $I(\CT ,\o)$ of non-isomorphic countable models
of a complete countable first-order theory $\CT$
behaves as the cardinality of ``nicely definable" (closed, Borel, analytic) subsets of Polish spaces:
$I(\CT ,\o)\leq \o$ or $I(\CT ,\o)= \c$.
The class of partial orders has a significant role in this context,
because by a result of Arnold Miller (see \cite{Stee}) the conjecture follows from its restriction to that class.
Regarding the class of partial orders we recall that Vaught's conjecture was confirmed for the theories of
linear orders (Matatyahu Rubin \cite{Rub}), model-theoretic trees (John Steel \cite{Stee}), reticles (James Schmerl \cite{Sch3}),
Boolean algebras (Paul Iverson \cite{Ive}).

Following the concept exploited in \cite{KMon,KAC,KFMD,Ksharp,KFLD}, we deal with the question
{\it whether Vaught's conjecture holds for (the theories of) the structures from $\la \CC\ra$,
if it holds for the structures from $\CC$},
where $\CC$ is a class of structures of a countable language
and $\la \CC\ra$ is a closure of $\CC$ obtained in some reasonable way.
For example,  Rubin in \cite{Rub} confirmed the ``sharp" version of Vaught's conjecture, VC$^\sharp$ (saying that $I(\CT,\o)\in \{1,\c\}$),
for linear orders with finitely many unary predicates,
while by \cite{KFMD} VC$^\sharp$ holds for all relational structures definable by quantifier free formulas in labeled linear orders
(i.e.\ for all structures from the closure $\la \CC^{\rm lo}_{\rm lab} \ra _{\rm def}$, called FMD structures\footnote{The structures from
$\la \CC^{\rm lo} \ra _{\rm def}$ are called {\it monomorphic} by Fra\"{\i}ss\'{e} \cite{Fra}
and $\la \CC^{\rm lo}_{\rm lab} \ra _{\rm def} $ is the class of structures admitting a {\it finite monomorphic decomposition} introduced by Pouzet and Thi\'{e}ry \cite{PT}.}).
In addition, by \cite{Ksharp} VC$^\sharp$ holds for all partial orders from the closure
of the class of linear orders $\CC^{\rm lo}$ under finite direct products and disjoint unions, $\la \CC^{\rm lo}\ra _{\du\Pi}$.

In \cite{Ksharp} VC$^\sharp$ was confirmed for the partial orders from the class $\la \CC\ra _{\du\Pi}$,
where $\CC$ is either the class of rooted FMD trees or the class of initially finite rooted trees satisfying VC$^\sharp$.
The key assumption providing these results was that deleting the root from a tree we obtain a tree having finitely many connected components,
which is a strong restriction leaving a large part of the class $\CC ^{\rm rt}$ of rooted trees unexplored
(e.g.\ one can ask what is going on with the product $\o ^{<\o}\times \o ^{<\o}$).

So, our initial intension was to confirm Vaught's conjecture
for the theories of partial orders belonging to the closure $\la \CC^{\rm rt}\ra _{\du\Pi}$
and, in particular, to understand what is going on with finite products of rooted trees.
But it turned out that much more can be said about the model theory of such structures,
which will be described in the sequel.

This paper deals with finite direct products $\prod _{i<n}\X_i$ and disjoint unions $\bcd _{i<n}\X_i$ of $L_b$-structures;
let us call these $\X _i$-s ``factors" in both cases and let $\square _{i<n}\X_i$ denote both constructions.
While these two constructions of models preserve the basic model-theoretic relations, $\cong$, $\equiv$ and $\preccurlyeq$,
in one direction (e.g., if the factors are elementarily equivalent, then the products are too),
it is natural to ask the following questions concerning the other direction.
\begin{itemize}
\item[Q1] {\it When $\square _{i<m}\Y_i \cong \square _{i<n}\X_i$ implies that $m=n$ and  $\X _i \cong \Y _{\pi (i)}$, for $i<n$, for some $\pi \in \Sym (n)$}?\\
          This sometimes fails: for Boolean algebras we have $P(1)\times P(1)\times P(2)\cong P(1)\times P(3)$.
\item[Q2] {\it When $\BE \preccurlyeq \square _{i<n}\X_i$ implies that $\BE=\square _{i<n}\BE_i$, where $\BE _i\preccurlyeq \X_i$, for $i<n$}?\\
          This sometimes fails too: If $\X =\la \o,\emptyset\ra$, then $\X ^2=\la \o\times \o, \emptyset\ra$
          and for each $E\in [\o\times \o]^\o$ we have $\BE :=\la E,\emptyset\ra \preccurlyeq \X ^2$.
          Or, the countable atomless Boolean algebra, $\B =\Clop (2^\o)$, elementarily embeds in its square, onto the diagonal
          (since $\Th (\B)$ is $\o$-categorical and model-complete).
\item[Q3] {\it When $\Y\equiv \square _{i<n}\X_i$ implies that $\Y \cong \square _{i<n}\Y_i$, where $\Y _i \equiv \X_i $, for $i<n$}? \\
          We note that replacing $\Y \cong \square _{i<n}\Y_i$ by $\Y = \square _{i<n}\Y_i$ we obtain a negative answer
          even when for $\square _{i<n}\X_i$ we take the square $\o \times \o$ of the ordinal $\o$ (which is a rooted tree); see Remark \ref{R002}.
\end{itemize}
In order to avoid repetition of similar proofs throughout the paper we consider
an ``abstract" operation $\square$ acting on the class $\Mod _L$ of models of a countable language $L$ and
a class $\CC \subset \Mod _L$ closed under $\equiv$
and in the following definition isolate conditions providing that, in particular,
Vaught's conjecture holds for the structures from the closure $\la \CC \ra_{\square}$ of $\CC$ under $\square$,
if it holds for the structures from $\CC$.
\begin{df}\label{D000}\rm
Let $L$ be a countable language
and let $\square $ be an operation
which for each $n\geq 1 $ to each $n$-tuple $\la \X_0, \dots ,\X _{n-1}\ra$ of $L$-structures
adjoins an $L$-structure $\square _{i<n}\X _i$
such that\footnote{Here we could be more pedantic
and regard an operation $\square$ on the class  $\Mod _L\!/\!\cong$ of isomorphism types, but this complicates the exposition and we avoid that approach.}
\begin{equation}\label{EQ621}\textstyle
\min \{\sum _{i<n} |X_i|,\prod _{i<n} |X_i| \}\leq |\square _{i<n}\X _i|\leq \max \{\sum _{i<n} |X_i|,\prod _{i<n} |X_i| \}.
\end{equation}
If $\CC$ is a class of $L$-structures closed under $\equiv$,
we will say that the pair $(\square ,\CC)$ is a {\it perfect pair} iff whenever $\X _i,\Y _i\in \CC$, for $i<\o$, and $m,n\geq 1$ we have
\begin{eqnarray}
\square _{i<n}\X _i \cong \square _{i<m}\Y _i & \Leftrightarrow & m=n \land \exists \pi \in \Sym (n)\;\forall i<n\;\X _i \cong \Y _{\pi (i)}, \label{EQ615}\\
\Y \equiv \square _{i<n}\X _i & \Leftrightarrow & \Y\cong \square _{i<n}\Y _i, \;\mbox{ where }\Y _i \equiv \X _i, \mbox{ for }i<n, \label{EQ616}\\
\BE \preccurlyeq \square _{i<n}\X _i & \Leftrightarrow & \BE = \square _{i<n}\BE _i, \;\mbox{ where }\BE _i \preccurlyeq \X _i, \mbox{ for }i<n. \label{EQ617}
\end{eqnarray}
Then we will write $\Perf (\square ,\CC)$.
We will write $\Perf _s (\square ,\CC)$ if, in addition, whenever $n\geq 1$, $\X _i\in \CC$, for $i<n$, and $\k$ is a cardinal we have
\begin{equation}\label{EQ618}
\forall i<n \;\; (\X _i \mbox{ is $\k$-saturated)} \;\Rightarrow\;  \square _{i<n}\X _i \mbox{ is $\k$-saturated.}
\end{equation}
\end{df}
\begin{rem}\label{R003}\rm
If $\Perf (\square ,\CC)$, $\X _i\in \CC$, for $i<n$, $\CT _i :=\Th (\X _i)$, for $i<n$, and $\CT:=\Th (\square _{i<n}\X _i)$, then
(\ref{EQ616}) provides a characterization of the models of $\CT$ via the models of $\CT _i$, $i<n$, namely,
\begin{equation}\label{EQ630}\textstyle
\Mod (\CT)=\bigcup \{ [\square _{i<n}\Y _i]: \forall i<n \;\Y _i \in \Mod (\CT _i)\},
\end{equation}
by (\ref{EQ615}) for the isomorphism relation on the class $\Mod (\CT)$ we have:
if $\Y\in [\square _{i<n}\Y _i]$ and $\Z \in [\square _{i<n}\Z _i]$, then
\begin{equation}\label{EQ631}
\Y \cong \Z \Leftrightarrow
\exists \pi \in \Sym (n)\;\forall i<n\;\Y _i \cong \Z _{\pi (i)}
\end{equation}
and by (\ref{EQ617}) for the preorder $\prec$ on $\Mod (\CT)/\!\cong$,
defined by $[\Y ]\prec [\Z]$ iff there is an elementary embedding $f:\Y \prec \Z$, we have
\begin{equation}\label{EQ632}
[\Y ]\prec [\Z]\Leftrightarrow \Y\cong \square _{i<n}\BE _i, \;\mbox{ where }\BE _i \preccurlyeq \Z _i, \mbox{ for }i<n.
\end{equation}
Thus the property $\Perf (\square ,\CC)$ provides simple descriptions of the model-theoretic phenomena in $\la \CC\ra _{\square}$.
\end{rem}
Moreover, in Section \ref{S3} we show that $\Perf (\square ,\CC)$ provides the following additional properties of $\CT$:
\begin{equation}\label{EQ633}\textstyle
\max \{I(\CT_i ):i<n \}\leq I(\CT )\leq \prod _{i<n}I(\CT_i ),
\end{equation}
so, if Vaught's conjecture is true for $\CT _i$, $i<n$, then it is true for $\CT $.
Further, $\CT$ is $\o$-categorical (resp.\ atomic) iff $\CT _i$, $i<n$, are $\o$-categorical (resp.\ atomic);
and, under $\Perf _s(\square ,\CC)$,
$\CT$ is small iff $\CT _i$, $i<n$, are small;
also, the countable atomic and the saturated model of $\CT$ (if they exist) have the expected form.

In Section \ref{S4} we prove first that $\Perf _s(\prod ,\CC ^{\rm rt}_{>1})$, where $\CC ^{\rm rt}_{>1}=\{ \X \in \CC ^{\rm rt}:|X|>1\}$,
and then show that the theory $\CT:=\Th (\prod _{i<n}\X_i)$ of a product of rooted trees $\prod _{i<n}\X_i \in \la \CC ^{\rm rt}\ra_{\Pi}$
has all nice properties listed above except (\ref{EQ615}) (since $|X_i|=1$ is allowed).
In particular, Vaught's conjecture for $\CT$ is confirmed.

In Section \ref{S5} we show that $\Perf _s(\bcd ,\CC ^{\rm fd})$, where $\CC ^{\rm fd}$ is the class of $L_b$-structures of finite diameter;
so we confirm that the theories $\Th (\bcd _{i<n}\X_i)$, where $\bcd _{i<n}\X_i \in \la \CC ^{\rm fd}\ra_{\du }$, have the aforementioned properties.

On the basis of the results obtained in Sections \ref{S4} and \ref{S5} in Section \ref{S6} we obtain analogous results
for the theories of the structures from the closure $\la \CC ^{\rm rt}\ra_{\du\Pi}$ of the class of rooted trees.
In particular, using Steel's result concerning theories of trees, Vaught's conjecture is confirmed for such theories.
\section{Preliminaries}\label{S2}
Throughout the paper we assume that $L$ is a countable language; in particular, let $L_b=\la R\ra$, where $R$ is a binary relational symbol.
$\Mod _L$ and $\Form _L$ will denote the class of $L$-structures and the set of all first-order $L$-formulas.
Isomorphism and elementary equivalence are denoted by $\cong$ and $\equiv$.
$I(\CT, \o):= |\Mod (\CT, \o)/\!\cong |$ is the number of non-isomorphic countable models of a complete theory $\CT$;
if $\CT$ has a finite model, $\X$, then, for convenience, we define $\Mod (\CT, \o)/\!\cong \;:= \{ [\X]\}$ and $I(\CT, \o):=1$.
For simplicity, instead of $I(\CT, \o)$ and $I(\Th (\X), \o)$ we will write only $I(\CT)$ and $I(\X )$.

If $\X,\Y \in \Mod _L$, then
$\X \preccurlyeq \Y $ denotes that $\X$ is an elementary substructure of $\Y$,
$f:\X \prec \Y$ denotes that $f$ is an elementary embedding
and $\X \prec \Y$ denotes that such mapping exists.
A model $\X$ is
{\it prime} iff $\X \prec \Y$, for each $\Y \equiv \X$;
{\it countably prime} iff $\X$ is countable and $\X \prec \Y$, for each countable $\Y \equiv \X$;
{\it atomic} iff any $n$-tuple $\bar x\in X$ satisfies a complete formula in $\Th (\X)$;
{\it $\k$-saturated} iff for each $A\in [X]^{<\k}$ the expansion $\X _A:=(\X ,\la a:a\in A\ra)$ realizes each 1-type of $L_A$ consistent with $\Th (\X _A)$;
{\it countably saturated} iff $\X$ is countable and $\o$-saturated;
{\it countably universal} iff $\X$ is countable and for each countable $\Y \equiv \X$ we have $\Y \prec \X$.
A complete $L$-theory $\CT$ is
{\it small} iff $|\bigcup _{n\in \o}S_n(\CT)|\leq \o$;
{\it atomic} iff each formula which is consistent with $\CT$ is completable in $\CT$.
We will use the following standard facts; for (a), (b) and (c) see \cite{CK}, pp.\ 98, 105, 94 and 96; claim (d) follows from definitions.
\begin{fac} \label{T538}
Let $|L|\leq \o$, let $\CT$ be a complete theory of $L$ and let $\X$ and $\Y$ be $L$-structures. Then

(a) $\CT$ is small iff $\CT$ has a countably saturated model iff $\CT$ has a countably universal model;

(b) $\CT$ has a countable atomic model iff $\CT$ is atomic;

(c) $\X$ is a countable atomic model iff $\X$ is a prime model iff $\X$ is a  countably prime model;

(d) If  $f:\X \rightarrow \Y$ is an isomorphism and $\BE \preccurlyeq \X$, then $f[\BE] \preccurlyeq \Y$.
\end{fac}
We recall basic facts concerning direct products (valid for any language $L$).
Claims (a) and (b) are evident; for (c) and (d) see \cite{FV} or \cite{Hodg}, p.\ 462;
for (e) see \cite{FV} or \cite{Hodg}, p.\ 486.
\begin{fac}\label{T043}
If $\X _i$ and $\Y _i$, for $i\in I$, are $L$-structures, then we have

(a) $\prod _{i\in I}\X _i \cong \prod _{i\in I}\X _{\pi (i)}$, for any permutation $\pi \in \Sym (I)$;

(b) If $\X _i \cong \Y _i$, for all $i\in I$, then $\prod _{i\in I}\X _i \cong \prod _{i\in I}\Y _i$;

(c) If $\X _i \equiv \Y _i$, for all $i\in I$, then $\prod _{i\in I}\X _i \equiv \prod _{i\in I}\Y _i$;

(d) If $\X _i \preccurlyeq \Y _i$, for all $i\in I$, then $\prod _{i\in I}\X _i \preccurlyeq \prod _{i\in I}\Y _i$;

(e) If $|I|<\o$ and $\X _i $, $i\in I$, are $\k$-saturated, then $\prod _{i\in I}\X _i$ is $\k$-saturated.
\end{fac}
If $\X =\la X,\r \ra$ is an $L_b$-structure,
then the transitive closure $\r _{\rm rst}$ of the relation $\r _{\rm rs}:=\Delta _X \cup \r \cup \r ^{-1}$,
given by $x \;\r_{\rm rst} \;y$ iff there are $n\in \N$ and $z_0 =x , z_1, \dots ,z_n =y$ such that $z_i \;\r _{\rm rs} \;z_{i+1}$, for each $i<n$;
(the tuple $\la z_0,\dots ,z_n\ra$ is a {\it path} from $x$ to $y$)
is the minimal equivalence relation on $X$ containing $\r$.
The equivalence classes $[x]_{\r _{\rm rst}}$, $x\in X$, are the {\it connected components} of $\X$
and $\X$ is said to be {\it connected} iff $|X/\r _{\rm rst}|=1$.
The formula $\f _{\rm rs}(u,v):= u=v \lor R(u,v) \lor R(v,u)$ defines the relation $\r_{\rm rs}$ in the structure $\X$
and for $n\in \N$ the formula
\begin{equation}\label{EQ581}\textstyle
\ve _n(u,v):=\exists w_0, \dots w_n \;( w_0=u \land w_n=v \land \bigwedge _{i<n} \f _{\rm rs}(w_i,w_{i+1}))
\end{equation}
says that there is a path of length $\leq n$ from $u$  to $v$.
$\X$ is said to be {\it of diameter $\leq n$}, we write $\delta (\X )\leq n$,
iff for each $x,y\in X$ there is a path of length $\leq n$ from $x$ to $y$,
iff $\X\models \f _{\delta \leq n}$, where $\f _{\delta \leq n}:=\forall u,v\;\ve _n(u,v)$.
{\it $\X$ is of finite diameter}, in notation $\delta (\X )<\o$, iff $\delta (\X )\leq n$, for some $n\in \N$.
It is easy to see that $\delta (\X )<\o$ iff all models of $\Th (\X )$ are connected (see Fact 2.11 of \cite{Ksharp}).

If $\X _i =\la X_i ,\r _i\ra$, $i\in I$, are  $L_b$-structures with pairwise disjoint domains,
then the $L_b$-structure $\bcd _{i\in I} \X _i :=\la \bigcup _{i\in I} X_i , \bigcup _{i\in I} \r_i \ra$
is called the {\it disjoint union} of the structures $\X _i$, $i\in I$.
\begin{fac}\label{T535}
If  $\X _i:i<n$ and $\Y _i:i<n$ are families of pairwise disjoint connected $L_b$-structures, then

(a) $\bcd _{i<n}\X _i \cong \bcd _{i<n}\Y _i $ iff there is $\pi\in \Sym (n)$ such that $\X_i \cong \Y _{\pi (i)}$, for each $i<n$;

(b) If $\X _i \equiv  \Y _i$, for all $i<n$, then $\bcd _{i<n}\X _i \equiv  \bcd _{i<n}\Y _i$;

(c) If $\X _i \preccurlyeq  \Y _i$, for all $i<n$, then $\bcd _{i<n}\X _i \preccurlyeq  \bcd _{i<n}\Y _i$;

(d) If $\X _i$, $i<n$, are $\k$-saturated, then $\bcd _{i<n}\X _i$ is $\k$-saturated.
\end{fac}
\dok
For (a) and (b) see \cite{Ksharp}, Fact 2.6 and Propositions 2.4(b).
For (c) see Fact 2.5(b) of \cite{KFLD}.

 (d) First we prove the claim for $n=2$.
Let $\X =\la X, R^{\X }\ra$ and $\Y =\la Y, R^{\Y }\ra$ be $\k$-saturated $L_b$-structures.
Let $L_0=\la S\ra$, $L_1=\la T\ra$ and $L_0+L_1 =\la S,T,P,Q\ra$,
where $S$, $T$, $P$ and $Q$ are new different relational symbols, $\ar (S)=\ar (T)=2$ and $\ar (P)=\ar (Q)=1$.
Then for $\X ':=\la X, S^{\X '}\ra$, where $S^{\X '}= R^{\X }$, and $\Y ':=\la Y, T^{\Y' }\ra$, where $T^{\Y '}= R^{\Y }$,
we have $\X '\in \Mod _{L_0}$, $\Y '\in \Mod _{L_1}$
and $\X '+\Y':=\la X \cup Y , S^{\X '}, T^{\Y' },X,Y\ra \in \Mod _{L_0+L_1}$ is the disjoint sum of  $\X '$ and $\Y'$ (see \cite{Hodg}, p.\ 101).
Since $\X '=\X$ and $\Y '=\Y$, the structures  $\X '$ and $\Y'$ are $\k$-saturated
and (see \cite{Hodg}, p.\ 488) $\X '+\Y'$ is a $\k$-saturated $(L_0+L_1)$-structure.
Let $\Gamma$ be the 1-dimensional unrelativized injective interpretation of the language $L_b$ in the language $L_0+L_1$ (see \cite{Hodg}, p.\ 216),
given by the $(L_0+L_1)$-formula $\f _R (u,v):=(P(u)\land P(v)\land S(u,v))\lor (Q(u)\land Q(v)\land T(u,v))$.
Then $\Gamma (\X '+\Y') =\X \;\dot{\cup }\;\Y $ is a $\k$-saturated $L_b$-structure (see \cite{Hodg}, p.\ 486).
The claim for $n$ structures is proved by induction
(note that the proof for $n=2$ is given for {\it any} two $L_b$-structures; we did not assume that they are connected).
\kdok
A partial order $\X =\la X,\leq \ra$ is a {\it (model-theoretic) tree}
iff its suborders $(\cdot ,x]:=\{ y\in X : y\leq  x\}$, $x\in X$, are linearly ordered.
If  $r\in X$ and $r\leq x$, for all $x\in X$,
then $r$ is the {\it root} of $\X$, $\X$ is a {\it rooted tree}
and $\X ^+:=\X \setminus \{ r \}$ is a tree (possibly disconnected).
If $\X=\la X, \leq \ra$ is a partial order and $r\not\in X$,
let $\X _r:=\la \{ r\} \cup X, \leq _r\ra$,
where $\leq _r \;=\{ \la r,x\ra :x\in \{ r\} \cup X\} \;\cup \leq$
(thus, $r=\min \X _r$).
\begin{fac}\label{T522}
If $\Y$ is a partial order, $r\not\in Y$ and $\X \preccurlyeq \Y$, then $\X _r\preccurlyeq \Y _r$.
\end{fac}
\dok
By Fact 2.5(b) of \cite{KFLD},
if $\sum _{\BI}\Y _i$ is a lexicographic sum of $L_b$-structures
and $\X _i \preccurlyeq  \Y _i$, for $i\in I$, then $\sum _{\BI}\X _i \preccurlyeq \sum _{\BI}\Y _i $.
Since $\Y _r=\sum _{\BI}\Y_i$, where $\BI =\la \{ 0,1\}, <\ra$, $\Y _0=\la \{ r\}, \{ \la r,r\ra\}\ra$ and $\Y_1 =\Y$,
and $\X _r=\sum _{\BI}\X_i$, where $\X _0=\Y _0$ and $\X_1 =\X$,
we have $\X _r\preccurlyeq \Y _r$.
\kdok
\paragraph{Definable sets and partitions}
The results of this paragraph will be used throughout the paper and are valid for any relational language $L$.
If $\X$ is an $L$-structure and $\f (\bar v)=\f (v_0,\dots ,v_{n-1})$ a first-order $L$-formula,
we use notation  $D_{\f (\bar v) ,\X}:=\{ \bar x \in X^n : \X \models \f [\bar x]\}$,
or $D_{\f ,\X}$, when the context admits. For $A\subset X$ by $\A$ we will denote the corresponding substructure of $\X$.
Sometimes we will introduce new constants
and use standard facts concerning that operation
(e.g.\ $\X\models \f[x_0,x_1]$ iff $(\X,x_0)\models \f (c,x_1)$).
\begin{fac}\label{T533}
If $\X$ and $\Y$ are $L$-structures and $\f (\bar v)\in \Form _L$, then

(a) If $\BE \preccurlyeq \X $, then $D_{\f ,\BE}= E^n\cap D_{\f ,\X}$;

(b) If $f:\X \rightarrow\Y$ is an isomorphism, then $f[D_{\f ,\X}]=D_{\f ,\Y}$.
\end{fac}
\dok
(a) Clearly, $\bar e \in D_{\f ,\BE}$
iff $\bar e\in E^n$ and $\BE\models \f [\bar e]$
iff (since $\BE \preccurlyeq \X$) $\bar e\in E^n$ and $\X\models \f [\bar e]$
iff $\bar e\in E^n \cap D_{\f ,\X}$.
Claim (b) is evident too.
\kdok
\begin{prop}\label{T519}
If $\X\in \Mod _L$ and $\ve (u,v) \in \Form _L$, where
$D_{\ve ,\X}$ is an equivalence relation on the set $X$
and $X/D_{\ve ,\X}=\{ X_i :i<n\}$, then

(a) If $\Y  \equiv \X$, then $D_{\ve ,\Y}$ is an equivalence relation on $Y$ and $|Y/D_{\ve ,\Y}|=n$;


(b) If $\BE\preccurlyeq \X$, then $E/D_{\ve ,\BE}=\{ E\cap X_i :i<n\}$ and $\BE\cap \X_i \preccurlyeq \X _i$, for all $i<n$;

(c) If $\Y \equiv \X$, then there is an enumeration $Y/D _{\ve ,\Y}=\{ Y_i :i<n\}$ such that $\Y _i\equiv \X_i$, for $i<n$.
\end{prop}
\dok
Statement (a) is a folklore fact (see e.g.\ \cite{KFMD}).

(b) By (a) $D_{\ve ,\BE}$ is an equivalence relation on the set $E$ and $|E/D_{\ve ,\BE}|=n$.
Let $C\in E/D_{\ve ,\BE}$, $e'\in C$ and $i<n$, where $e'\in X_i$.
Then we have
$e \in C$
iff $e\in E$ and $\la e,e'\ra \in D_{\ve ,\BE}$
iff $e\in E$ and $\BE\models \ve [e,e']$
iff $e\in E$ and $\X \models \ve [e,e']$ (since $\BE\preccurlyeq \X$)
iff $e\in E$ and $e,e' \in X_j$, for some $j<n$ (since $X/D_{\ve ,\X}=\{ X_i :i<n\}$)
iff  $e\in E\cap X_i$ (since $e'\in X_i$).
So, $C =E\cap X_i$ and the inclusion $E/D_{\ve ,\BE}\subset\{ E\cap X_i :i<n\}$ is proved.
Since $|E/D_{\ve ,\BE}|=n$ we have the equality.
Consequently, for each $i<n$ we have $E\cap X_i\neq \emptyset$.

Let $e\in E\cap X_i$
and let $L_c:=\la L , c\ra$, where $c$ is a new constant.
Since $\BE \preccurlyeq \X$ and $e\in E $ we have $(\BE ,e) \preccurlyeq _{L_c} (\X ,e)$.
For the $L_{c}$-formula $\f (v):=\ve (c,v)$ and $x\in X$ we have
$x\in D_{\f , (\BE ,e)}$,
iff $x\in E$ and $(\BE ,e)\models \f [x]$,
iff (since  $(\BE ,e) \preccurlyeq _{L_c} (\X ,e)$) $x\in E$ and $(\X ,e)\models \ve (c,v) [x]$,
iff $x\in E$ and $\X \models \ve [e,x]$,
iff $x\in E $ and $\la e,x\ra \in D_{\ve ,\X}$,
iff (since $e\in X_i$) $x\in E $ and $x\in X_i$.
Thus
\begin{equation}\label{EQ578}
D_{\f , (\BE ,e)}=E \cap X_i\;\;\mbox{ and }\;\; D_{\f , (\X ,e)}= X_i ,
\end{equation}
where the second equality is obtained as above, taking $\X$ instead of $\BE$.

Let $U$ be a new unary relational symbol and $L^+ :=\la L_c , U\ra$.
By (\ref{EQ578}) the $L^+$-structure $(\BE ,e, E\cap X_i)$ is a definitional expansion of $(\BE ,e)$
and the $L^+$-structure $(\X ,e, X_i)$ is a definitional expansion of $(\X ,e)$.
So, since $(\BE ,e) \preccurlyeq _{L_c} (\X ,e)$ we have (see \cite{Hodg} p.\ 60)
$\BE ^+:=(\BE ,e, E\cap X_i) \preccurlyeq _{L^+} (\X ,e, X_i) =:\X^+$,
which for the corresponding relativised (to $U$) reducts (to $L$)
gives $\BE ^+ _U \preccurlyeq \X ^+ _U$ (see \cite{Hodg} p.\ 203),
where $\BE ^+ _U$ is the substructure of the reduct $\BE ^+|L =\BE \subseteq \X$ with the domain $U^{\BE ^+}=E\cap X_i$; that is, $\BE ^+ _U=\BE \cap \X_i$,
while $\X ^+ _U$ is the substructure of the reduct $\X ^+|L =\X$ with the domain $U^{\X ^+}=X_i$; that is $\X ^+ _U=\X_i$.
So, $\BE \cap \X_i\preccurlyeq \X_i$ indeed.

(c) Let $\max\{|X|,|Y|\}=\k \geq \o$ and let $\Z$ be a $\k ^+$-saturated model such that $\X,\Y \preccurlyeq\Z$.
By (a) $D_{\ve ,\Z}$ is an equivalence relation on the set $Z$ and $|Z/D_{\ve ,\Z}|=n$;
also, $D_{\ve ,\Y}$ is an equivalence relation on the set $Y$ and $|Y/D_{\ve ,\Y}|=n$.
By (b) (applied to $\Z$ as the larger structure),
for $i<n$ there is an element of $Z/D_{\ve ,\Z}$, say $C_i$, such that $X_i=X\cap C_i$ and $\X _i \preccurlyeq \C _i$,
and there is an element of $Y/D_{\ve ,\Y}$, say $Y_i$, such that $Y_i=Y\cap C_i$ and $\Y _i \preccurlyeq \C _i$;
thus  $\Y _i\equiv \X_i$, for $i<n$.
In addition, if $i\neq j$, then $X_i\neq X_j$ and, hence, $C_i\neq C_j$,
which gives $Y_i\neq Y_j$; so $Y/D_{\ve ,\Y}=\{ Y_i:i<n\}$.
\kdok
\section{Perfect pairs}\label{S3}
In this section we assume that $|L|\leq \o$, that $\square $ is an operation
which for each $n\geq 1 $ to each $n$-tuple $\la \X_0, \dots ,\X _{n-1}\ra \in \Mod _L^n$ adjoins a structure $\square _{i<n}\X _i \in \Mod _L$
and that $\CC \subset \Mod _L$ is a class closed under $\equiv$.
Recall that the properties $\Perf (\square ,\CC)$ and $\Perf _s(\square ,\CC)$ were established in Definition \ref{D000}.
\begin{te}\label{T543}
Let $\Perf (\square ,\CC)$, let $\X :=\square  _{i<n}\X_i$, where $\X _i\in \CC$, for $i<n$, let $\CT:=\Th (\X )$ and $\CT _i:=\Th (\X _i)$, for $i<n$.
Then \\[-6mm]
\begin{itemize}\itemsep -0.4mm
\item[\rm (a)] $\max \{I(\CT_i ):i<n \}\leq I(\CT )\leq \prod _{i<n}I(\CT_i ) $;
               thus, $I(\CT )=1$  iff $\,I(\CT_i)=1$, for each $i<n$;
\item[\rm (b)] If  the theories $\CT_i$, $i<n$, satisfy Vaught's conjecture, then $\CT$ satisfies Vaught's conjecture: defining $\k:=\prod _{i<n}I(\CT _i)$ we have
               $I(\CT )=\k$, if $\k\in \{1,\o,\c \}$, and $I(\CT ) \in [3,\o )$, otherwise;
\item[\rm (c)] $\CT$ is atomic iff $\;\CT _i$, $i<n$, are atomic;
               then $\square _{i<n}\X_i^{\rm at}$ is a countable atomic model of $\CT$, where $\X_i^{\rm at}$ is a countable atomic model of $\CT _i$, for $i<n$;
\item[\rm (d)] If, in addition, $\Perf _s(\square ,\CC)$, then $\CT$ is small iff $\;\CT _i$, $i<n$, are small;
               then $\square  _{i<n}\X_i^{\rm sat}$ is a countably saturated model of $\CT$, where $\X _i ^{\rm sat}$ is a countably saturated model of $\CT_i$, for $i<n$.
\end{itemize}
\end{te}
\dok
(a) Let $\sim$ be the equivalence relation on the set $n$ defined by
\begin{equation}\label{EQ611}
i\sim j \Leftrightarrow \CT _i = \CT_j.
\end{equation}
Let $j<n$ and let us choose $\Y _i \in \Mod (\CT_i,\o)$, for $i\in n\setminus [j]$.
For $\A \in \Mod (\CT _j ,\o)$ let $\Z ^\A =\square  _{i<n}\Z_i^\A$, where for $i<n$ we define
\begin{equation}\label{EQ612}\textstyle
           \Z _i^\A=\left\{
           \begin{array}{cl}
           \A, & \mbox{ if }\;i\in [j];\\
           \Y_i, & \mbox{ if }\;i\in n\setminus[j].
           \end{array}
              \right.
\end{equation}
Then $\Z_i^\A \equiv \X _i$, for $i<n$,
by (\ref{EQ616}) we have $\Z^\A \equiv \X $
and by (\ref{EQ621})  $\Z^\A\in \Mod (\CT  ,\o)$.
Let $\B \in \Mod (\CT _j ,\o)$, where $\B \not\cong \A$;
then $\Z^\B\in \Mod (\CT  ,\o)$ too.
Suppose that $\Z^\A\cong\Z^\B$.
By (\ref{EQ615}) there is $\pi \in \Sym (n)$ such that $\A =\Z _j^\A \cong \Z _{\pi(j)}^\B$,
which implies that $\pi (j)\in n\setminus[j]$ and, hence, $\Z _{\pi(j)}^\B =\Y_{\pi(j)}$.
Thus $\X_j \equiv\A \equiv \Y_{\pi(j)}\equiv \X_{\pi(j)}$, which is false because $j\not\sim \pi (j)$.
So, if $\A,\B \in \Mod (\CT _j ,\o)$ and $\A\not\cong \B$, then $\Z^\A,\Z^\B\in \Mod (\CT  ,\o)$ and $\Z^\A\not\cong\Z^\B$.
Consequently, $I(\CT )\geq I(\CT _j)$, for all $j<n$, and the first inequality is proved.

Concerning the second inequality, if $|X|<\o$, then $I(\CT)=1$ and we are done.
Let $|X|\geq \o$, let $\l _i:=\min \{|X_i|,\o\}$ and $\k_i:=I(\CT_i )$,  for $i<n$,
and let $\Mod (\CT_i,\o)/\!\cong \;=\{ [\X _i^\a]:\a < \k _i\}$ be an enumeration.
We prove that $\Mod (\CT,\o)/\!\cong \;=\{ [\square  _{i<n}\X _i^{\a _i}]: \la \a_i :i<n\ra \in \prod _{i<n}\k _i\}$,
which will give $I(\CT )\leq \prod _{i<n}I(\CT_i)$.
So, if $\Y \in \Mod (\CT,\o)$,
then by (\ref{EQ616}) we have $\Y \cong \square  _{i<n}\Y _i$,
where for each $i<n$ we have $\Y _i \equiv \X _i$;
by (\ref{EQ621}) we have $|Y_i|\leq \o$ and, hence, $|Y_i|=\l _i$;
so, $\Y _i \cong \X _i^{\a _i}$, for some $\a _i <\k _i$,
and by (\ref{EQ615}) we have $\Y \cong \square  _{i<n}\X _i^{\a _i}$, where $\la \a_i :i<n\ra \in \prod _{i<n}\k _i$.
Conversely, if $\la \a_i :i<n\ra \in \prod _{i<n}\k _i$,
then  $\X _i^{\a _i}\equiv \X _i$, for $i<n$,
by (\ref{EQ616}) we have $\square  _{i<n}\X _i^{\a _i}\equiv \square  _{i<n}\X _i =\X$,
and by (\ref{EQ621}) $|\square  _{i<n}\X _i^{\a _i}|\leq \o$.
Thus, $[\square  _{i<n}\X _i^{\a _i}]\in \Mod (\CT,\o)/\!\cong$.

(b) By the assumption we have $I(\CT)\in [1,\o]\cup \{ \c\}$, for all $i<n$.
By (a), if $I(\CT_i )=\c$, for some $i<n$, then $I(\CT )=\c$;
otherwise we have $I(\CT_i)\leq \o$, for all $i<n$.
Now, if $I(\CT_i )=\o$, for some $i<n$, then by (a) we have $I(\CT )=\o$;
and, by (a) again,  $I(\CT )=1$ iff $I(\CT_i)=1$, for all $i<n$.
The remaining case is when $I(\CT_i)=n_i\in \N$, for all $i<n$,
and (by Vaught's theorem) $n_i\geq 3$, for some $i<n$; then by (a) $3\leq I(\CT )<\o $.

(c) Let $\CT $ be an atomic theory and let $\X ^{\rm at}$ be a countable atomic model of $\CT$.
By (\ref{EQ616}) we have $\X ^{\rm at}\cong \square _{i<n}\Y_i$, where $\Y _i \equiv \X_i$, for all $i<n$,
and w.l.o.g.\ we assume that $\X ^{\rm at}= \square _{i<n}\Y_i$.
By (\ref{EQ621}) we have $|Y_i|\leq \o$, for $i<n$.
Let $\sim$ be the equivalence relation on the set $n$ defined by (\ref{EQ611}).
Let $j\in n$ and $\A \in \Mod (\CT _j ,\o)$
and let $\Z ^\A =\square  _{i<n}\Z_i^\A$, where the factors $\Z _i^\A $, $i\in n$, are defined by (\ref{EQ612}).
Then $\Z_i^\A \equiv \Y _i$, for $i<n$,
by (\ref{EQ616}) we have $\Z^\A \equiv \X ^{\rm at}$
and, by (\ref{EQ621}), $\Z ^\A \in \Mod (\CT  ,\o)$.
By Fact \ref{T538}(c) $\X ^{\rm at}$ is a prime model of $\CT$
so, there is $f: \X ^{\rm at} \prec \Z^\A $
and, hence, $f[\X ^{\rm at}]\preccurlyeq \Z^\A$.
By (\ref{EQ617}) we have $f[\X ^{\rm at}]=\square  _{i<n}\BE_i$, where $\BE _i \preccurlyeq \Z^\A_i$, for $i<n$,
which implies that $\square  _{i<n}\Y_i\cong\square  _{i<n}\BE_i$
and by (\ref{EQ615}) there is $\pi \in \Sym (n)$ such that
$\Y_i \cong \BE _{\pi (i)} \preccurlyeq  \Z^\A _{\pi (i)}$, for each $i<n$.
So, if $i\in \pi ^{-1} [[j]]$, that is, $\pi (i)\in [j]$,
then by (\ref{EQ612}) we have $\Y _i\prec \Z^\A _{\pi (i)}=\A$,
which gives $i\sim j$, that is $i\in [j]$.
Thus $\pi ^{-1} [[j]]\subset [j]$
and, since $[j]$ is a finite set, $\pi^{-1} [[j]]= [j]$.
So, since $j\in \pi^{-1} [[j]]$, we have  $\Y_j \prec  \A$.
Thus $\Y _j \prec  \A$, for each $\A \in \Mod (\CT _j ,\o)$,
which means that $\Y_j$ is a countably prime model of $\CT_j$,
by Fact \ref{T538}(c) $\Y_j$ is a countable atomic model of $\CT_j$
and by Fact \ref{T538}(b) the theory $\CT_j$ is atomic.

Conversely, if the theories $\CT _i$, $i<n$, are atomic,
then by Fact \ref{T538}(c) for each $i<n$ there is a prime model of $\CT _i$, say $\X_i^{\rm at}$.
If $\Y \equiv \X$,
then by (\ref{EQ616}) we have $\Y \cong \square  _{i<n}\Y_i$, where $\Y _i\equiv \X_i$, for $i<n$.
So, $\X_i^{\rm at}\prec \Y _i$, for $i<n$,
by (\ref{EQ617}) we have $\square _{i<n}\X_i^{\rm at}\prec \square  _{i<n}\Y_i$
and, hence, $\square _{i<n}\X_i^{\rm at}$ is a prime model of $\CT$.
By Fact \ref{T538}(c) $\square _{i<n}\X_i^{\rm at}$ is a countable atomic model of $\CT$
and by Fact \ref{T538}(b) the theory $\CT$ is atomic.

(d) Let $\Perf _s(\square ,\CC)$ and let $\CT$ be a small theory;
by Fact \ref{T538}(a) there is a countably saturated  model of $\CT$, say $\X ^{\rm sat}$.
Then $\X ^{\rm sat}\equiv \X$,
and by (\ref{EQ616}) $\X ^{\rm sat}\cong\square  _{i<n}\Y_i$, where $\Y _i\equiv \X_i$, for $i<n$;
w.l.o.g.\ suppose that $\X ^{\rm sat}=\square  _{i<n}\Y_i$.
By (\ref{EQ621}) we have $|Y_i|\leq \o$, for $i<n$.
Let $\sim$ be the equivalence relation on the set $n$ defined by (\ref{EQ611}).
Let $j\in n$ and $\A \in \Mod (\CT _j ,\o)$
and let $\Z ^\A =\square  _{i<n}\Z_i^\A$, where the structures $\Z _i^\A$, for $i\in n$, are defined by (\ref{EQ612}).
Then $\Z_i^\A \equiv \Y _i$, for $i<n$,
by (\ref{EQ616}) we have $\Z^\A \equiv \X ^{\rm sat}$,
so, there is $f: \Z^\A \prec \X ^{\rm sat}$
and, hence, $f[\Z^\A]\preccurlyeq \X ^{\rm sat}$.
By (\ref{EQ617}) we have $f[\Z^\A]=\square  _{i<n}\BE_i$, where $\BE _i \preccurlyeq \Y_i$, for $i<n$,
which implies that $\square  _{i<n}\Z_i^\A\cong\square  _{i<n}\BE_i$
and by (\ref{EQ615}) there is $\pi \in \Sym (n)$ such that
$\Z_i^\A \cong \BE _{\pi (i)} \preccurlyeq  \Y _{\pi (i)}$, for each $i<n$.
So, for $i\in [j]$ by (\ref{EQ612}) we have $\A =\Z_i^\A \prec  \Y _{\pi (i)}$,
which gives $j\sim \pi (i)$, that is $\pi (i)\in [j]$.
Thus $\pi [[j]]\subset [j]$
and, since $[j]$ is a finite set, $\pi [[j]]= [j]$.
So, since $j\in [j]$, we have $j=\pi (i)$, for some $i\in [j]$
and $\A \prec  \Y _{\pi (i)}=\Y _j$.
Thus $\A \prec  \Y _j$, for each $\A \in \Mod (\CT _j ,\o)$,
which means that $\Y_j$ is a countably universal model of $\CT_j$
and by Fact \ref{T538}(a) the theory $\CT_j$ is small.

Conversely, if the theories $\CT _i$, $i<n$, are small,
by Fact \ref{T538}(a) there are countably saturated models $\X _i ^{\rm sat}\models\CT_i$, $i<n$,
by (\ref{EQ621}) and (\ref{EQ618}) the structure $\square  _{i<n}\X_i^{\rm sat}$ is countably saturated,
and, since  $\X _i ^{\rm sat}\equiv \X_i$, for $i<n$, by (\ref{EQ616}) we have $\square  _{i<n}\X_i^{\rm sat}\equiv \square  _{i<n}\X_i$,
that is, $\square  _{i<n}\X_i^{\rm sat}\models \CT$.
\kdok
\begin{prop}\label{T546}
If (\ref{EQ621}), (\ref{EQ615}), (\ref{EQ617}) and (\ref{EQ618}) hold and
whenever $n\geq 1$ and $\X _i,\Y _i \in \CC$, for $i<n$,  we have
\begin{equation}\label{EQ622}
(\forall i<n \;\X _i \equiv \Y _i) \Rightarrow \square _{i<n}\X _i \equiv \square _{i<n}\Y _i ,
\end{equation}
\noindent
then (\ref{EQ616}) is true and hence we have $\Perf _s(\square ,\CC)$.
\end{prop}
\dok
The implication ``$\Leftarrow$" in (\ref{EQ616}) follows from (\ref{EQ622}) and we prove ``$\Rightarrow$".
If $|X_i|<\o$, for all $i<n$,
then by (\ref{EQ621}) the structure $\square _{i<n} \X _i$ is finite,
$\Y \equiv\square _{i<n} \X _i$ implies that $\Y\cong \square _{i<n} \X _i$ and we are done.
So, suppose that $|X_i|\geq\o$, for some $i<n$, and let $\Y \equiv \square _{i<n} \X _i$; then by (\ref{EQ621}) $|Y|\geq\o$.
If $\k:=\max\{|Y|,|X_0|,\dots,|X_{n-1}|\}$,
then for each $i<n$ there is a $\k ^+$-saturated model $\X _i^s \succcurlyeq \X _i$
and, since $\X _i^s\in \CC$, for $i<n$, by (\ref{EQ618}) the structure $\square _{i<n} \X _i^s$ is $\k ^+$-saturated as well.
Thus $\X _i^s \equiv \X _i$, for $i<n$,
by (\ref{EQ622}) we have $\square _{i<n} \X _i^s \equiv \square _{i<n} \X _i$
and, hence, $\Y \equiv\square _{i<n} \X _i^s$,
which, by the $\k ^+$-universality of $\square _{i<n} \X _i^s$ (see \cite{Hodg}, p.\ 482), gives $\Y \preccurlyeq\square _{i<n} \X _i^s$.
By (\ref{EQ617}) there are $\Y _i \preccurlyeq \X _i^s$ such that $\Y =\square _{i<n} \Y _i$
and, clearly, $\Y _i \equiv \X _i$, for $i<n$.
\kdok
\begin{rem}\label{R001}\rm
By Theorem \ref{T543}(a) the equality $I(\CT )=\prod _{i<n} I(\CT _i )$ holds, if $\prod _{i<n} I(\CT _i )\geq \o$, or if $\prod _{i<n} I(\CT _i )=1 $.
In order to prove that, in addition,
$$\textstyle
( \CT _i \neq \CT _j \mbox{ for different }i,j<n) \Rightarrow I(\CT )=\prod _{i<n} I(\CT _i ),
$$
we show that the enumeration $\Mod (\CT,\o)/\!\cong \;=\{ [\square  _{i<n}\X _i^{\a _i}]: \bar \a \in \prod _{i<n}\k _i\}$
given in the proof of Theorem \ref{T543}(a) is one-to-one.
So, if $\bar \a ,\bar \b \in \prod _{i<n}\k _i$ and $\square _{i<n} \X _i^{\a _i}\cong \square _{i<n} \X _i^{\b _i}$,
then by (\ref{EQ615}) there is $\pi \in \Sym (n)$ such that for each $i<n$ we have $\X _i^{\a _i}\cong \X _{\pi (i)}^{\b _{\pi (i)}}$,
which implies $\CT _i=\CT _{\pi (i)}$ and, hence, $\pi (i)=i$.
Thus $\pi =\id _n$ and for each $i<n$ we have $\X _i^{\a _i}\cong \X _i^{\b _i}$,
which gives $\a _i =\b _i$;
so, $\bar \a =\bar \b$.
\end{rem}
\begin{rem}\rm
In Definition \ref{D000} we do not request that $\CC$ is closed under $\square$;
in fact, our intention is to consider the closure $\la \CC\ra_{\square}$ of $\CC$ under $\square$.
In order to provide that $\CC \subset\la \CC\ra_{\square}$ we could request that $\square _{i<1}\X _0 =\X _0$
and in order to provide that the closure is obtained in one step,
roughly, that $\la \CC\ra_{\square}=\{\square _{i\in I}\X _i: |I|<\o \land \la \X_i:i\in I\ra \in \CC ^I \}$,
we should define $\square _{i\in I}\X _i$ for finite $I$
and then request that acting on $\Mod _L$ ``a $\square $-sum of $\square $-sums is a $\square $-sum";
namely, that $\square _{i\in I}\square _{j\in J_i}\X _{i,j} \cong \square _{\la i,j\ra\in \bigcup _{i\in I} \{ i\}\times  J_i}\X _{\la i,j\ra}$.
But all the statements in this section are true without these assumptions (which are satisfied by finite products and disjoint unions considered in the sequel).
\end{rem}
\section{Products of rooted trees}\label{S4}
For a proof that Vaught's conjecture holds for products of rooted trees
we intend to apply Theorem \ref{T543} but there is an obstacle:
the pair $(\prod ,\CC ^{\rm rt})$ is not
perfect (see Definition \ref{D000})
because (\ref{EQ615}) fails if in a product $\prod _{i<n}\X_i$ we allow factors of size 1.
In this section we will first show that $\Perf _s(\prod ,\CC ^{\rm rt}_{>1})$
and apply Theorem \ref{T543} to that pair;
then we will extend that result to the class $\la\CC ^{\rm rt}\ra _{\Pi}$.
\subsection{The closure $\la\CC ^{\rm rt}_{>1}\ra _{\Pi}$: isomorphism and unique factorization}
Here we show that in the class $\la\CC ^{\rm rt}_{>1}\ra _{\Pi}$ the number of factors is a first-order property
and that the factorization is unique up to isomorphism and a permutation of factors (which is expected by Fact \ref{T043}(a)).
\begin{te}\label{T505}
If $\X _i$, $i<n$, and $\Y _j$, $j<m$, are rooted trees of size $>1$, where $m,n\in \N$, then we have
\begin{itemize}\itemsep -0.4mm
\item[\rm (a)] There are $L_b$-sentences $\l _k$, $k\in \N$, such that  $\prod _{i<n}\X_i\models \bigwedge _{1\leq k\leq n}\l_k \land \neg \l _{n+1}$;
\item[\rm (b)] $\prod _{i<n}\X_i \equiv \prod _{j<m}\Y_j \Rightarrow m=n$;
\item[\rm (c)] $\prod _{i<n}\X_i \cong \prod _{j<m}\Y_j \Leftrightarrow m=n \land \exists \pi \in \Sym (n)\;\;\forall i<n\;\; \X _i \cong \Y _{\pi (i)} $.

In addition, if $f:\prod _{i<n}\X_i \rightarrow \prod _{i<n}\Y_i$ is an isomorphism
and $I_\pi : \prod _{i<n}\Y_i \rightarrow \prod _{i<n}\Y_{\pi (i)}$ is the isomorphism defined by $I_\pi (\la y_0,\dots ,y_{n-1}\ra)=\la y_{\pi (0)},\dots ,y_{\pi (n-1)}\ra$,
then  there are isomorphisms $f_i :\X _i \rightarrow \Y _{\pi (i)}$, $i<n$,
such that $I_\pi\circ f=\prod _{i<n} f_i$.
\end{itemize}
\end{te}
A proof of the theorem is given at the end of this section. First we construct a convenient decomposition of products of rooted trees.

Let $\X :=\prod _{i<n}\X_i$, where $n\in \N$ and $\X _i$, $i<n$, are rooted trees of size $>1$.
For $i<n$, let $r_i$ be the root of $\X _i$,
$X_i^+ :=X_i\setminus \{ r_i\}$
and let $\{X _i^j :j\in J_i\}$ be the partition of the substructure $\X _i^+$ of $\X _i$ into its connectivity components.
Then $\X _i ^+ =\bcd _{j\in J_i}\X _i^j$ is a disjoint union of connected trees.
Defining $X^{(m)}=\{ \bar x\in \prod _{i<n}X_i: |\{ i<n:x_i>r_i \}|=m\}$, for $m\leq n$, we have
\begin{equation}\label{EQ589}\textstyle
X^{(m)}=\bigcup _{K\in[n]^m} \bigcup _{\bar j\in \prod _{i\in K}J_i}\Big\{ \bar x\in \prod _{i<n}X_i : \forall i\in K \;(x_i\in X_i^{j_i}) \land \forall i\in n\setminus K \;(x_i=r_i)\Big\},
\end{equation}
we obtain a partition $\{ X^{(m)}:0\leq m\leq n\}$ of the set $X$
and we will show that it is definable in $\X$.
In particular, $X^{(0)}=\{ \bar r\}$, where $\bar r=\la r_0,\dots ,r_{n-1}\ra$,
and $X^{(0)}=D_{\f _0,\X}$, where $\f _0 (v):= \forall w \;\;v\leq w$.
For $m=1$ by (\ref{EQ589}) we have

\begin{eqnarray}
X^{(1)} & = &\textstyle\bigcup _{i<n}\bigcup _{j\in J_i}S_i^j,
              \;\;\mbox{ where } S_i^j :=\{ r_0\} \times \dots \times \{ r_{i-1}\} \times X _i^j\times \{ r_{i+1}\} \times \dots \times \{ r_{n-1}\},\label{EQ590}\\
        & = &\textstyle
             \bigcup _{i<n}A_i,\;\mbox{ where } A _i:= \bigcup _{j\in J_i} S_i^j =\{ r_0\} \times \dots \times \{ r_{i-1}\} \times X _i^+\times \{ r_{i+1}\} \times \dots \times \{ r_{n-1}\}. \label{EQ623}
\end{eqnarray}
\begin{cla}\label{T507}\rm
$X ^{(1)}= D_{\f _1,\X }$, where $\f _1 (v)$ is a formula saying that $v$ is not a smallest element and that the set $(\cdot ,v]$ is linearly ordered, say
\begin{equation}\label{EQ592}\textstyle
\f _1 (v) :=  \exists u \;(u<v)\land  \forall u,w \;(u,w \leq v \Rightarrow u\leq w \lor w\leq u).
\end{equation}
\end{cla}
\dok
If $\bar x\in X ^{(1)}$, then by (\ref{EQ590}) there are $i<n$ and $j\in J_{i}$
such that $\bar x =\la r_0, \dots, r_{i-1}, x_i,  r_{i+1}, \dots , r_{n-1}\ra $, where $r_{i}<x_i\in  X _i^j$,
and, hence, $\bar r <_\X \bar x$.
If $\bar y,\bar z \leq_{\X } \bar x$,
then for $i'\in n \setminus \{ i \}$ we have $y_{i'}, z_{i'} \leq_{\X _{i'} } x_{i'}=r_{i'}$
and, hence, $y_{i'}= z_{i'} =r_{i'}$.
Since $y_i, z_i\leq_{\X _{i}} x_i$ and $\X _i$ is a tree, we have:
either $y_i\leq _{\X _i} z_i$, which gives $\bar y\leq _{\X }\bar z$;
or $z_i\leq _{\X _i} y_i$, which gives $\bar z\leq _{\X }\bar y$.
Thus $\X \models \f _1[\bar x]$.

If $\bar x\not\in X ^{(1)}$,
then $\bar x\in S$, where $\dim (S)\neq 1$.
So, either $\bar x =\bar r$
and, hence,  $\X \not\models \f _1[\bar x]$,
or there are different $i',i'' <n$
such that $r_{i'}<_{\X _{i'}}x_{i'}$ and $r_{i''}<_{\X _{i''}}x_{i''}$.
Then for $\bar y=\la r_0, \dots, r_{i'-1}, x_i',  r_{i'+1}, \dots , r_{n-1}\ra$
and $\bar z=\la r_0, \dots, r_{i''-1}, x_i'',  r_{i''+1}, \dots , r_{n-1}\ra$
we have $\bar y,\bar z\leq_{\X } \bar x$,
but $\bar y\not\leq_{\X } \bar z$ and $\bar z\not\leq_{\X } \bar y$;
so $\X \not\models \f _1[\bar x]$ again.
\kdok
\begin{cla}\label{T532}\rm
For the following substructures of the product $\X$  we have
\begin{equation}\label{EQ591}\textstyle
\S_i^j \cong  \X _i^j , \mbox{ for }j\in J_i, \;\;\mbox{ and }\;\;
\A _i = \bcd _{j\in J_i} \S_i^j \cong \X _i ^+  \;\;\mbox{ and }\;\;
\X^{(1)}=\bcd _{i<n}\A_i.
\end{equation}
\end{cla}
\dok
By (\ref{EQ590}) the mapping $f:\X _i^j \rightarrow \S _i^j$
defined by $f(x)=\la r_0, \dots,r_{i-1}, x,r_{i+1},\dots,r_{n-1}\ra$
is an isomorphism and the first claim is true.
For the second, since $\X _i ^+=\bcd _{j\in J_i}\X _i^j$ is a disjoint union of trees,
it remains to be proved that for different $j,j'\in J_i$, $\bar x\in S _i^j$ and $\bar y\in S _i^{j'}$
we have $\bar x \not\leq _\X \bar y$.
But then $\bar x=\la r_0, \dots,r_{i-1}, x,r_{i+1},\dots,r_{n-1}\ra$, where $x\in X _i^j$,
and $\bar y=\la r_0, \dots,r_{i-1}, y,r_{i+1},\dots,r_{n-1}\ra$, where $y\in X _i^{j'}$,
and $\bar x \leq _\X \bar y$ would imply that $x\leq _{\X _i ^+}y$, which is false.
So, the second claim is true.
By (\ref{EQ590}) we have $X^{(1)}=\bcd _{i<n}A_i$
so it remains to show that for different $i,i'<n$, $\bar x\in S _i^j$ and $\bar y\in S _{i'}^{j'}$
we have $\bar x \not\leq _\X \bar y$.
But then $\bar x=\la r_0, \dots,r_{i-1}, x,r_{i+1},\dots,r_{n-1}\ra$, where  $x\in X _i^j$,
and $\bar y=\la r_0, \dots,r_{i'-1}, y,r_{i'+1},\dots,r_{n-1}\ra$, where  $y\in X _{i'}^{j'}$,
so $\bar x \leq _\X \bar y$ would imply that $x=r_i$, which is false.
\kdok
\begin{cla}\label{T508}\rm
If $\r(u,v)$ is a formula
saying that $u$ and $v$ are in $X^{(1)}$ and have no upper bound in $\X$ or are compatible in $X^{(1)}$, say
\begin{equation}\label{EQ601}\textstyle
\r (u,v):= \f _1(u) \land \f _1(v) \land \Big(\neg \exists w \;\; (u\leq w \land v\leq w) \lor \exists w \;(\f _1 (w)\land w\leq u \land w\leq v) \Big),
\end{equation}
then $D _{\r , \X }= \bigcup _{i<n} A_i ^2 $ (that is, $D _{\r , \X }$ is the equivalence relation on the set $X^{(1)}$ determined by the partition $\{ A_i:i<n\}$).
\end{cla}
\dok
Let $\la \bar x, \bar y \ra\in X^2$ and $\X \models \r [\bar x, \bar y]$.
Then by Claim \ref{T507} we have $\bar x, \bar y \in X^{(1)}$
and by (\ref{EQ591}) $\bar x\in A_i$ and $\bar y\in A_{i'}$, for some $i,i' <n$.
Suppose that $i\neq i'$, say $i<i'$.
Then there are $j\in J_i$, $x \in X_i^j$, $j'\in J_{i'}$ and  $y\in X_{i'}^{j'}$ such that
\begin{equation}\label{EQ567}
\bar x=\la r_0, \dots,r_{i-1}, x,r_{i+1},\dots,r_{n-1}\ra \;\;\mbox{ and }\;\;\bar y=\la r_0, \dots,r_{i'-1}, y,r_{i'+1},\dots,r_{n-1}\ra.
\end{equation}
So for $\bar z:=\la r_0, \dots,r_{i-1}, x,r_{i+1},\dots, r_{i'-1}, y,r_{i'+1},\dots ,r_{n-1}\ra$ we have $\bar x ,\bar y \leq _\X \bar z$
and, since $\X \models \r [\bar x, \bar y]$,
there is $\bar u \in X^{(1)}$ such that $\bar u \leq _\X \bar x,\bar y$.
But by (\ref{EQ567}) $\bar u \leq _\X \bar x,\bar y$ implies that $\bar u=\bar r \not\in X^{(1)}$, and we have a contradiction.
Thus $i'=i$ and $\la \bar x, \bar y \ra\in A_i^2$.

Conversely, if $i<n$ and $\la \bar x, \bar y \ra\in A_i^2$,
then there are $j,j'\in J_i$, $x\in X_i^j$ and $y\in X_i^{j'}$ such that
\begin{equation}\label{EQ595}
\bar x=\la r_0, \dots,r_{i-1}, x,r_{i+1},\dots,r_{n-1}\ra \;\;\mbox{ and }\bar y=\la r_0, \dots,r_{i-1}, y,r_{i+1},\dots,r_{n-1}\ra
\end{equation}
and by (\ref{EQ590}) and  Claim \ref{T507} we have  $\X \models \f _1[\bar x]$ and $\X \models \f _1[\bar y]$.
If there exists $\bar z \in X$ such that $\bar x, \bar y \leq _\X \bar z$,
then by (\ref{EQ595}) $x,y \leq _{\X _i}z_i$,
which implies that $j=j'$ and that $x$ and $y$ are comparable in $\X _i$.
So, if $x \leq _{\X _i}y$, then $\bar x \leq _\X \bar y$
and, hence, $\bar x$ and $\bar y$ are compatible in $\X^{(1)}$; so, $\X \models \r [\bar x, \bar y]$ .
The same holds if $y \leq _{\X _i}x$.
\kdok
Let $\a _m (v_0,\dots ,v_m)$ be a formula saying that $v_m=\sup \{v_0, \dots ,v_{m-1}\}$, say
$\a _m := \bigwedge _{k < m}v_k\leq v_m \land \forall u \;(\bigwedge _{k< m}v_k\leq u \Rightarrow v_m \leq u)$
and for $m\in [2,n]$ let $\f _m (v)$ be a formula saying that $v$ is a supremum of $m$ elements of $X^{(1)}$ belonging to different sets $A_i$, say
\begin{equation}\label{EQ604}\textstyle
\f _m (v)  :=  \exists v_0,\dots ,v_{m-1} \;\Big( \bigwedge _{k<m} \f _1(v_k)\land \bigwedge _{k<l<m} \neg \r (v_k,v_l) \land \a _m (v_0,\dots ,v_{m-1},v)\Big).
\end{equation}
\begin{cla}\label{T528}\rm
$X^{(m)}=D_{\f _m ,\X}$, for $2\leq m\leq n$.
\end{cla}
\dok
If $\bar x\in X^{(m)}$,
then by (\ref{EQ589}) there is $K=\{ i_0 ,\dots ,i_{m-1}\}\in [n]^m$ such that
$x_i =r_i$, for $i\in n\setminus K$,
and for $k<m$ there is $j_k\in J_{i_k}$ such that $x_{i_k}\in X_{i_k}^{j_k}$;
thus, $\bar y ^k :=\la r_0, \dots ,r_{i_k-1}, x_{i_k}, r_{i_k+1}, \dots, r_{n-1}\ra \in A_{i_k}$.
So by Claims  \ref{T507} and \ref{T508} we have $\X \models \f _1 [\bar y^k]$, for $k<m$, and $\X \models \neg \r [\bar y^k,\bar y^l]$, for different $k,l<m$.
It is evident that $\bar y ^k \leq _\X \bar x$, for all $k<m$.
If $\bar z \in X$ and $\bar y ^k \leq _\X \bar z$, for all $k<m$,
then for $i\in n\setminus K$ we have $x_i =r_i \leq _{\X _i} z_i$,
and for $k<m$ we have $x_{i_k}=y_{i_k}\leq _{\X _{i_k}} z_{i_k}$;
so, $\bar x \leq _\X \bar z$
and, hence, $\bar x =\sup _\X \{\bar y^k : k<m\}$.
Thus $\X \models \a _m[\bar y ^0,\dots ,\bar y ^{m-1},\bar x]$ and $\X \models \f _m[\bar x]$, that is $\bar x\in D_{\f _m ,\X}$.
Conversely, if $\X \models \f _{m}[\bar x]$,
then there are $\bar y^0,\dots ,\bar y^{m-1}\in X^{(1)}$ and different $i_k \in n$ such that $\bar y^k\in A_{i_k}$, for $k<m$,
and $\bar x =\sup \{\bar y^0,\dots, \bar y^{m-1}\}$.
For $k<m$ we have $\bar y ^k :=\la r_0, \dots ,r_{i_k-1}, y_{i_k}, r_{i_k+1}, \dots, r_{n-1}\ra $, where $y_{i_k}\in X_{i_k}^{j_k}$, for some $j_k\in J_{i_k}$,
and defining $z_i=r_i$, for $i\in n\setminus \{i_k:k<m\}$, and $z_{i_k}=y_{i_k}$, for $k<m$,
we have $\bar z\in  X^{(m)}$ and, as above, $\bar z =\sup \{\bar y^0,\dots, \bar y^{m-1}\}$.
Consequently, $\bar x =\bar z\in  X^{(m)}$.
\kdok
\begin{cla}\label{T534}\rm
$\CP_\X := \{ X^{(m)}: m \in \{ 0\} \cup \{2,\dots ,n\}\} \cup \{ A_i:i<n\}$ is a partition of the set $X$
and $\CP_\X =X/D_{\ve ,\X}$, where $D_{\ve ,\X}$ is the equivalence relation on the set $X$ defined in $\X$ by the formula
\begin{equation}\label{EQ602}\textstyle
\ve (u,v):=\bigvee _{m \in \{ 0\} \cup \{2,\dots ,n\}}(\f_m (u)\land \f_m (v)) \lor \r (u,v).
\end{equation}
\end{cla}
\dok
By (\ref{EQ589}) and (\ref{EQ590}) $\CP_\X $ is a partition of the set $X$
and the corresponding equivalence relation, say $\sim$ is given by $\bar x \sim \bar y$ iff $\bar x ,\bar y \in P$, for some $P\in \CP_\X$.
By (\ref{EQ602})and Claims \ref{T507}, \ref{T508} and \ref{T528} for $\bar x ,\bar y \in X$ we have
$\X\models \ve [\bar x ,\bar y]$ iff $\bar x ,\bar y\in X^{(m)}$, for some $m \in \{ 0\} \cup \{2,\dots ,n\}$,
or $\bar x ,\bar y\in X^{(1)}$ and $\bar x ,\bar y\in A_i$, for some $i<n$. Thus $D_{\ve ,\X} =\;\sim $.
\kdok

\paragraph{Proof of Theorem \ref{T505}}
Let $\X _i$, $i<n$, and $\Y _j$, $j<m$, be rooted trees of size $>1$, where $m,n\in \N$.

(a) For $n\in \N$ let $\l _k :=\exists v_0,\dots ,v_{k-1} \;( \bigwedge _{i<k} \f _1(v_i)\land \bigwedge _{i<j<k} \neg \r (v_i,v_j))$.
Since for each $i<n$ we have $|X_i|>1$  we can take $\bar x ^i\in A_i$, for $i<n$.
by Claims \ref{T532} and \ref{T507} we have $\X \models \f _1[\bar x^i]$, for $i<n$,
and, by Claim \ref{T508} $\X \models \neg \r [\bar x^i,\bar x^j ]$, whenever $i<j<n$;
thus $\X \models \l _n$. Similarly, $\X \models \l _k$, whenever $1\leq k \leq n$.
Assuming that $\X \models \l _{n+1}$
by the same claims there would be $\bar x^i \in X^{(1)}$, for $i\leq n$,
such that $\X \models \neg \r [\bar x^i,\bar x^j ]$, whenever $i<j\leq n$,
which is false by Claim \ref{T508};
thus $\X \models \neg \l _{n+1}$.

(b) We prove the contrapositive. If $m\neq n$, say $m<n$,
then by (a)  we have  $\prod _{i<n}\X_i\models \l _{m+1}$ and  $\prod _{j<m}\Y_j\models \neg \l_{m+1}$;
thus $\prod _{i<n}\X_i \not\equiv \prod _{j<m}\Y_j$.

(c) The implication ``$\Leftarrow$" follows from (a) and (b) of Fact \ref{T043}.
For the converse we suppose that $\X :=\prod _{i<n}\X_i \cong \prod _{j<m}\Y_j =:\Y $,
note that $m=n$ follows from (b)
and take an isomorphism $f:\X \rightarrow \Y$.
By Claim \ref{T534} the formula $\ve (u,v)$ given by (\ref{EQ602})
in the posets $\X$ and $\Y$ defines equivalence relations $D_{\ve ,\X}$ and $D_{\ve ,\Y}$ on the sets $X$ and $Y$ respectively
corresponding to their partitions
\begin{eqnarray*}
\CP_\X := X/D_{\ve ,\X} & = & \{ X^{(m)}: m \in \{ 0\} \cup \{2,\dots ,n\}\} \cup \{ A_i^{\X}:i<n\},\\
\CP_\Y := Y/D_{\ve ,\Y} & = & \{ Y^{(m)}\,: m \in \{ 0\} \cup \{2,\dots ,n\}\} \cup \{ A_i^{\Y}:i<n\} .
\end{eqnarray*}
By Fact \ref{T533}(b) we have $f[D_{\ve ,\X}]=D_{\ve ,\Y}$,
and, hence, for each set $P\in\CP_\X$ we have $f[P]\in \CP_\Y$.
Moreover, for $m\in \{0,\dots, n\}$ by Claim \ref{T507}, \ref{T528}  and Fact \ref{T533}(b) we have $f[X^{(m)}]=f[D_{\f _m,\X}]=D_{\f _m,\Y}=Y^{(m)}$.
In particular, $f[X^{(1)}]=Y^{(1)}$,
which implies that $f$ maps $A_i^\X$-s onto $A_i^\Y$-s,
namely there is $\pi\in \Sym (n)$ such that $f[A_i^{\X}]=A_{\pi (i)}^{\Y}$, for $i<n$.
Thus for each $i<n$ we have $\A_i^{\X} \cong \A_{\pi (i)}^{\Y}$
and, by Claim \ref{T532},  $\X_i^+ \cong \Y _{\pi (i)} ^+$,
which clearly implies $\X_i \cong \Y _{\pi (i)} $.

Further, let $I_\pi : \prod _{i<n}\Y_i \rightarrow \prod _{i<n}\Y_{\pi (i)}$,
where $I_\pi (\bar y)=\la y_{\pi (0)},\dots ,y_{\pi (n-1)}\ra$, for $\bar y \in Y$.
Clearly $f(\bar r ^\X)=\bar r ^\Y$, where $\bar r ^\X :=\la r^{\X _i}:i<n\ra$ and $\bar r ^\Y :=\la r^{\Y _i}:i<n\ra$;
thus, $f[\{ \bar r ^\X\} \cup A_i^{\X}]=\{ \bar r ^\Y\}\cup A_{\pi (i)}^{\Y}$, for $i<n$.
Let $i<n$.
If $x\in X_i$,
then by (\ref{EQ623}) $\bar x^i:=\la r^{\X _0}, \dots ,r^{\X _{i-1}}, x, r^{\X _{i+1}},\dots , r^{\X _{n-1}}\ra \in \{ \bar r ^\X\} \cup A_i^{\X}$
so, $f(\bar x^i)\in \{ \bar r ^\Y\}\cup A_{\pi (i)}^{\Y}$,
which gives $f(\bar x^i)=\la r^{\Y _0}, \dots ,r^{\Y _{\pi(i)-1}}, y, r^{\Y _{\pi (i)+1}},\dots , r^{\Y _{n-1}}\ra$,
and, hence, $p_{\pi (i)}(f(\bar x^i))=y\in Y_{\pi (i)}$, where $p_{\pi (i)}:\prod _{i<n}Y_i\rightarrow Y_{\pi (i)}$ is a projection.
So, the mapping $f_i :X_i \rightarrow Y_{\pi (i)}$ given by
\begin{equation}\label{EQ627}
f_i (x)=p_{\pi (i)}(f(\la r^{\X _0}, \dots ,r^{\X _{i-1}}, x, r^{\X _{i+1}},\dots , r^{\X _{n-1}}\ra)), \quad \mbox{ for }x\in X_i,
\end{equation}
is a bijection.
Since $x\leq _{\X _i} x'$
iff $\bar x^i \leq _\X \bar x'^i$
iff $f(\bar x^i) \leq_\Y f(\bar x'^i)$
iff $p_{\pi (i)}(f(\bar x^i)) \leq_{\Y _{\pi (i)}} p_{\pi (i)}(f(\bar x'^i))$,
$f_i :\X_i \rightarrow \Y_{\pi (i)}$ is an isomorphism.
Let $\bar x\in X$ and let $\bar x^i:=\la r^{\X _0}, \dots ,r^{\X _{i-1}}, x_i, r^{\X _{i+1}},\dots , r^{\X _{n-1}}\ra $, for $i<n$.
Then by (\ref{EQ623}) for $i<n$ we have $\bar x^i \in \{ \bar r ^\X\} \cup A_i^{\X}$ and, hence, $f(\bar x^i)\in \{ \bar r ^\Y\}\cup A_{\pi (i)}^{\Y}$,
which gives $f(\bar x^i)=\la r^{\Y _0}, \dots ,r^{\Y _{\pi (i)-1}}, y, r^{\Y _{\pi (i)+1}},\dots , r^{\Y _{n-1}}\ra $, where $y\in Y_{\pi (i)}$,
and by (\ref{EQ627}) $f_i (x_i)=p_{\pi (i)}(f(\bar x^i))=y$.
Thus $f(\bar x^i)=\la r^{\Y _0}, \dots ,r^{\Y _{\pi (i)-1}},f_i (x_i) , r^{\Y _{\pi (i)+1}},\dots , r^{\Y _{n-1}}\ra $.
It is easy to check that $\bar x=\bigvee _{i<n}\bar x^i$,
which gives $f(\bar x)=\bigvee _{i<n}f(\bar x^i)=\bigvee _{i<n}\la r^{\Y _0}, \dots ,r^{\Y _{\pi (i)-1}}, f_i (x_i), r^{\Y _{\pi (i)+1}},\dots , r^{\Y _{n-1}}\ra$.
So the $i$-th coordinate of $I_{\pi} (f(\bar x))$ is the $\pi (i)$-th coordinate of $f(\bar x)$, that is $f_i (x_i)$,
which is the $i$-th coordinate of $(\prod _{i<n} f_i)(\bar x)$. Thus $I_{\pi} (f(\bar x))=(\prod _{i<n} f_i)(\bar x)$, for all $\bar x\in X$.
\kdok
\subsection{The closure $\la\CC ^{\rm rt}_{>1}\ra _{\Pi}$: elementary substructures}
\begin{te}\label{T511}
If $\X _i$, $i<n$, are rooted trees of size $>1$, then for each $L_b$-structure $\BE$ we have
\begin{equation} \label{EQ600}\textstyle
\BE \preccurlyeq \prod _{i<n}\X_i \;\Leftrightarrow\; \BE= \prod_{i<n}\BE_i, \;\mbox{ where }  \BE _i \preccurlyeq \X _i, \mbox{ for }i<n.
\end{equation}
\end{te}
\dok
Let $\BE \preccurlyeq \X :=\prod _{i<n}\X_i$.
By Claim \ref{T534} we have
$X/D_{\ve ,\X} =\{ X^{(m)}:m\in \{0\} \cup \{2,\dots, n\}\}\cup \{ A_i:i<n\}$,
by Proposition \ref{T519}(b) $E/D_{\ve ,\BE} =\{ E\cap X^{(m)}:m\in \{0\} \cup \{2,\dots, n\}\}\cup \{E\cap A_i:i<n\}$ and
\begin{equation} \label{EQ598}
\BE\cap \X^{(m)}\preccurlyeq  \X^{(m)}, \mbox{ for all } m\leq n, \;\;\mbox{ and }\;\; \BE\cap \A _i\preccurlyeq  \A _i, \mbox{ for all } i<n.
\end{equation}
By Claim \ref{T507} we have $X^{(1)}=D_{\f _1, \X}$ and by Fact \ref{T533}(a) $D_{\f _1,\BE}=E \cap D_{\f _1, \X}=E \cap X^{(1)}$; thus
\begin{equation} \label{EQ605}
E \cap X^{(1)}=D_{\f _1,\BE}.
\end{equation}
For $m\geq 1$ let $\t_m$ be a sentence saying that  for each $m$ elements of $X^{(1)}$ belonging to different sets $A_i$ there exists a supremum, say
\begin{equation}\label{EQ603}\textstyle
\t _m  := \forall v_0,\dots ,v_{m-1}\;\Big( \bigwedge _{k<m} \f _1(v_k)\land \bigwedge _{k<l<m} \neg \r (v_k,v_l) \Rightarrow \exists  v\;\; \a _m (v_0,\dots ,v_{m-1},v)\Big).
\end{equation}
\begin{cla}\label{T524}\rm
$\X \models \t _m$, for $1\leq m \leq n$.
\end{cla}
\dok
Let $\bar x^0,\dots ,\bar x^{m-1}\in X^{(1)}$
and $\X \models \neg \r [\bar x^k,\bar x^l]$, for different $k,l<m$.
By Claims \ref{T507} and \ref{T508} there are different $i_k <n$, for $k<m$, such that $\bar x^k \in A_{i_k}$,
and, hence, $\bar x^k =\la r_0, \dots ,r_{i_{k}-1}, x_{i_k}, r_{i_k +1}, \dots, r_{n-1}\ra$, where $x_{i_k}\in X_{i_k}^j$, for some $j\in J_{i_k}$.
Let $\bar y\in X$, where $y_{i_k} = x_{i_k}$, for $k<m$, and $y_i=r_i$, for $i\in n \setminus \{ i_k :k<m\}$.
Then, clearly, $\bar x^k \leq _\X \bar y$, for all $k<m$.
If $\bar z \in X$ and $\bar x^k \leq _\X \bar z$, for all $k<m$, then for $k<m$ we have $z_{i_k} \geq _{\X _{i_k}} x_{i_k} =y_{i_k}$
and for $i\in n \setminus \{ i_k :k<m\}$ we have $z_i \geq _{\X _i} r_i =y_i$ as well;
so, $\bar y \leq _\X \bar z$.
Thus, $\bar y =\sup _\X \{ \bar x^0,\dots ,\bar x^{m-1} \}$, that is, $\X \models \a _m [\bar x^0,\dots ,\bar x^{m-1},\bar y]$.
\kdok
\begin{cla}\label{T525}\rm
If $\BE '\preccurlyeq \X $ and $E \cap X^{(1)}=E' \cap X^{(1)}$, then $E =E'$.
\end{cla}
\dok
Let $E \cap X^{(1)}=E' \cap X^{(1)}$.
Since $\{ X^{(m)}: m\leq n\}$ is a partition of the set $X$
on order to prove that $E=E'$ we show that $E \cap X^{(m)}=E' \cap X^{(m)}$, for all $m\leq n$.
First $X^{(0)}=\{ \bar r\}$; so, by (\ref{EQ598}) we have $E \cap X^{(0)}=E' \cap X^{(0)}=\{ \bar r\}$.
So, by Claim \ref{T528} it remains to be shown that for each $m\in \{2, \dots, n\}$ we have
\begin{equation}\label{EQ570}
E \cap D_{\f _m ,\X}=E' \cap D_{\f _m ,\X}.
\end{equation}
If $e\in E \cap D_{\f _m ,\X}$,
Then $\X \models \f _m [e]$
and, since $\BE\preccurlyeq \X $, $\BE \models \f _m [e]$;
so, there are $e_0,\dots ,e_{m-1}\in E$ such that,
first, for each $k<m$ we have $\BE\models \f _1 [e_k]$
that is $e_k\in D_{\f _1,\BE}$
and, by (\ref{EQ605}), $e_k\in \bigcup _{i<n}A_i$;
so, $e_k \in A_{i_k}$, for some $i_k <n$.
Second, for $k<l<m$ we have $\BE\models  \neg \r [e_k,e_l]$,
which gives $\X\models  \neg \r [e_k,e_l]$; thus, by Claim \ref{T508}, $i_k \neq i_l$.
Third, $\BE\models \a _m [e_0,\dots ,e_{m-1},e]$ and, since $\BE\preccurlyeq \X $, $\X\models \a _m [e_0,\dots ,e_{m-1},e]$; so, we have
\begin{equation}\label{EQ568}
e_0,\dots ,e_{m-1}\leq _\X e \;\;\mbox{ and }\;\;\forall x\in X \;(e_0,\dots ,e_{m-1}\leq _\X x \Rightarrow e \leq_\X x).
\end{equation}
Since $E \cap X^{(1)}=E' \cap X^{(1)}$ we have $e_k \in E'\cap A_{i_k}$, for $k<m$.
For $k<m$ we have $\X\models \f _1[e_k]$ and, hence, $\BE'\models \f _1[e_k]$
and for $k<l<m$ we have $\X\models  \neg \r [e_k,e_l]$, which gives $\BE'\models  \neg \r [e_k,e_l]$.
By Claim \ref{T524} we have $\BE' \models \t _m$
and, hence, there is $e'\in E'$ such that $\BE '\models  \a _m [e_0,\dots ,e_{m-1},e']$ and, hence, $\X\models  \a _m [e_0,\dots ,e_{m-1},e']$; that is,
\begin{equation}\label{EQ569}
e_0,\dots ,e_{m-1}\leq _\X e' \;\;\mbox{ and }\;\;\forall x\in X \;(e_0,\dots ,e_{m-1}\leq _\X x \Rightarrow e' \leq_\X x).
\end{equation}
Now, by (\ref{EQ568}) and (\ref{EQ569}) we obtain $e \leq_\X e'$ and $e' \leq_\X e$;
thus $e=e'\in E'$ and we have proved that $E \cap D_{\f _m ,\X}\subset E' \cap D_{\f _m ,\X}$.
In the same way we prove another inclusion.
\kdok
For $i<n$ let $\pi _i :X \rightarrow X_i$ be the projection given by $\pi _i (\bar x)=x_i$, for all $\bar x\in X$.
For $i<n$
let $(\A _i )_{\bar r}$ and $(\BE \cap \A _i)_{\bar r}$ be the substructures of $\X$ with domains
$$
(A _i )_{\bar r}:=\{ \bar r\} \cup A_i \;\;\mbox{ and }\;\; (E \cap A _i)_{\bar r}:=\{ \bar r\} \cup (E \cap A _i).
$$
\begin{cla}\label{T518'}\rm
For each $i<n$ we have $E_i :=  \pi _i [(E \cap A_i)_{\bar r}]\subset X_i$ and

(a) The surjective restriction  $p_i:=\pi _i \mid (A_i)_{\bar r}  :(\A _i)_{\bar r} \rightarrow \X_i$ is an isomorphism;

(b) $p_i[(E \cap A_i)_{\bar r}]=E_i$;

(c) $\BE _i \preccurlyeq \X_i$;

(d) $\BE = \prod _{i<n}\BE _i$.
\end{cla}
\dok
(a)
The mapping $p_i$ is a bijection
because $(A_i)_{\bar r}  =  \{ r_0\}\times \dots \times \{r_{i-1}\} \times X _i\times \{r_{i+1}\} \times \dots \times \{r_{n-1}\} $.
If $\bar x,\bar x' \in (A_i)_{\bar r}$,
then $\bar x \leq _\X\bar x'$
iff $x_i \leq _{\X _i}x_i'$
iff $p_i (\bar x ) \leq _{\X _i}p_i(\bar x')$;
so $p_i$ is an isomorphism.

(b)  Since $p_i:=\pi _i \mid (A_i)_{\bar r}$
and $(E \cap A_i)_{\bar r}\subset (A_i)_{\bar r}$
we have  $p_i[(E \cap A_i)_{\bar r}]=\pi_i[(E \cap A_i)_{\bar r}]=E_i$.

(c) By (\ref{EQ598}) we have $\BE \cap \A _i \preccurlyeq \A _i$,
by Fact \ref{T522} $(\BE \cap \A _i )_r\preccurlyeq (\A _i)_r $
and, by (a) and (b) and Fact \ref{T538}(d), $\BE _i =p_i[(E \cap A_i)_{\bar r}]\preccurlyeq p_i[(A_i)_{\bar r}] =\X_i$.

(d) By (c) and Fact \ref{T043}(d) we have $\prod _{i<n}\BE _i\preccurlyeq \X $.
So, by Claim \ref{T525} it remains to be proved that $(\prod _{i<n}E _i) \cap X^{(1)} =E \cap X^{(1)}$, that is
$$\textstyle
(\prod _{i<n}E _i) \cap \bigcup _{i<n}A_i =E \cap \bigcup _{i<n}A_i.
$$
Let $i_0<n$ and $\bar x\in (\prod _{i<n}E _i) \cap A_{i_0}$.
Since $\bar x\in A_{i_0}$ we have $\bar x=\la r_0, \dots, r_{i_0-1}, x ,r_{i_0+1}, \dots,r_{n-1} \ra$, where $x\in X_{i_0}^+$,
and, since $\bar x\in (\prod _{i<n}E _i) $ we have $x\in E _{i_0}$;
so, by the definition of $E _{i_0}$,
there is $\bar y= \la r_0, \dots, r_{i_0-1},y,r_{i_0+1}, \dots,r_{n-1} \ra \in E \cap A_{i_0}$ such that $x=\pi _{i_0}(\bar y)=y$
and, hence, $\bar x= \bar y\in E \cap A_{i_0}\subset E \cap \bigcup _{i<n}A_i$.
Conversely, if $i_0<n$ and $\bar x\in E  \cap A_{i_0}$,
then, since $\bar x\in A_{i_0}$ we have $\bar x=\la r_0, \dots, r_{i_0-1},x,r_{i_0+1}, \dots,r_{n-1} \ra$, where $x\in X_{i_0}^+$,
and, since $\bar x\in E \cap A_{i_0}$, we have $x=p_{i_0}(\bar x)\in E_{i_0}$.
By (b) we have $r_i \in E_i$, for all $i\in n\setminus \{ i_0\}$;
so,  $\bar x\in  \prod _{i<n}E _i $
and, hence, $\bar x\in (\prod _{i<n}E _i) \cap \bigcup _{i<n}A_i$.
\kdok
By Claim \ref{T518'} the implication ``$\Rightarrow$" in (\ref{EQ600}) is true.
The implication ``$\Leftarrow$" is Fact \ref{T043}(d).
\kdok
\begin{te}\label{T549}
$\Perf_s(\prod ,\CC^{\rm rt}_{>1})$; so, (a)--(d) of Theorem \ref{T543} are true for each $\prod  _{i<n}\X_i\in\la \CC^{\rm rt}_{>1}\ra_{\Pi}$.
\end{te}
\dok
The claim follows from Proposition \ref{T546}:
the class $\CC ^{\rm rt}_{>1}$ is first-order definable and, hence, closed under $\equiv$;
if $\X _i\in \CC ^{\rm rt}_{>1}$, for  $i<n$,
then (\ref{EQ621}) is true and
by Theorems \ref{T505}(c) and \ref{T511} and Fact \ref{T043}(e)
the operation $\prod$ acting on the class $\CC ^{\rm rt}_{>1}$ satisfies conditions (\ref{EQ615}), (\ref{EQ617}) and (\ref{EQ618}),
while (\ref{EQ622}) is Fact \ref{T043}(c).
\kdok
\subsection{The closure $\la \CC ^{\rm rt}\ra_\Pi$}
\begin{te}\label{T544}
If $\X =\prod  _{i<n}\X_i$, where $\X _i\in\CC ^{\rm rt}$, for $i<n$, $\CT:=\Th (\X )$ and $\CT _i:=\Th (\X _i)$, $i<n$, then\\[-5mm]
\begin{itemize}\itemsep -0.5mm
\item[\rm (a)] $\BE \preccurlyeq \prod _{i<n} \X _i$ iff $\;\BE = \prod _{i<n}\BE _i$, where $\BE _i\preccurlyeq \X _i$, for each $i<n$;
\item[\rm (b)] $\Y \equiv \prod _{i<n} \X _i $ iff $\;\Y \cong\prod _{i<n} \Y _i$, where $\Y _i \equiv \X _i$, for each $i<n$.
\item[\rm (c)] The theory $\CT$ satisfies Vaught's conjecture: $I(\CT )=\k:=\prod _{i<n}I(\CT _i)$, if $\k\in \{1,\o,\c \}$, and $I(\CT ) \in [3,\o )$, otherwise;
\item[\rm (d)] $\CT$ is atomic iff $\;\CT _i$, $i<n$, are atomic;
               then $\prod _{i<n}\X_i^{\rm at}$ is a countable atomic model of $\CT$, where $\X_i^{\rm at}$ is a countable atomic model of $\CT _i$, for $i<n$;
\item[\rm (e)] $\CT$ is small iff $\;\CT _i$, $i<n$, are small;
               then $\prod  _{i<n}\X_i^{\rm sat}$ is a countably saturated model of $\CT$, where $\X _i ^{\rm sat}$ is a countably saturated model of $\CT_i$, for $i<n$.
\end{itemize}
\end{te}
\dok
If $|X_i|>1$, for all $i<n$, then all statements follow from Theorems \ref{T549}, \ref{T543}
and Steel's result (Vaught's conjecture for trees);
if $|X|=1$, all statements are evident. So, in the sequel we assume that
\begin{equation}\label{EQ634}
\emptyset \neq J:=\{ i<n : |X_i|>1\}=m<n \;\;\mbox{ and }\;\;X_i=\{ x_i\}, \mbox{ for }i\in [m,n).
\end{equation}

(a) Let $\BE \preccurlyeq \X$.
The mapping $\pi _m:\X \rightarrow \prod _{i<m}\X _i $ defined by $\pi _m (\bar x)=\la x_0,\dots ,x_{m-1}\ra$ is an isomorphism,
by Fact \ref{T538}(d) we have $\pi_m [\BE]\preccurlyeq \prod _{i<m}\X _i$
and by Theorem \ref{T511} we have $\pi_m [\BE] =\prod _{i<m}\BE _i$, where $\BE _i\preccurlyeq \X_i$, for $i<m$.
Defining $\BE _i:=\X _i$,  for $i\in [m,n)$, by Fact \ref{T043}(d) we have $\BE ':=\prod _{i<n}\BE _i \preccurlyeq \prod _{i<n}\X _i$
and $\pi _m [\BE ']=\prod _{i<m}\BE _i=\pi_m [\BE]$,
which since $\pi _m$ is a bijection gives $\BE =\BE '$.
Thus, $\BE = \prod _{i<n}\BE _i$, where $\BE _i\preccurlyeq \X _i$, for each $i<n$.
The converse is Fact \ref{T043}(d).

(b) Let $\Y \equiv \X$.
By (\ref{EQ634}) we have $\X \cong \prod _{i<m}\X _i =:\Z\in \la \CC ^{\rm rt}_{>1}\ra _{\Pi}$;
so $\Y \equiv \Z$
and by Theorem \ref{T549} $\Y \cong\prod _{i<m} \Y _i$, where $\Y _i \equiv \X _i$, for $i<m$.
Defining $\Y _i :=\X _i$, for $i\in [m,n)$ we have $\Y \cong\prod _{i<n} \Y _i$, where $\Y _i \equiv \X _i$, for $i<n$.
The converse follows from Fact \ref{T043}(c).

(c), (d) and (e).
By (\ref{EQ634}) we have $\X \cong \prod _{i<m}\X _i =:\Y$,
thus $\CT=\Th (\Y)$,
and for $i\in [m,n)$ we have: $\CT_i$ is atomic, small and $I(\CT _i)=1$.
Thus $\prod _{i<n}I(\CT_i )=\prod _{i<m}I(\CT_i )$
and, since (c) is true for $\Th(\Y)$,
for $\k \in \{ 1,\o,\c\}$ we have $I(\CT)=\k$ iff $I(\Th (\Y))=\k$ iff  $\prod _{i<m}I(\CT_i )=\k$ iff $\prod _{i<n}I(\CT_i )=\k$;
thus  (c) is true for $\CT$.
Also, $\CT$ is atomic iff $\Th (\Y)$ is atomic iff $\CT_i$, $i<m$, are atomic iff $\CT_i$, $i<n$, are atomic;
so (d) is true for $\CT$ and (e) has a similar proof.
\kdok
\section{The closure $\la \CC ^{\rm fd}\ra _{\du}$ (disjoint unions of $L_b$-structures of finite diameter)}\label{S5}
Here we show that $\Perf _s (\bcd ,\CC ^{\rm fd})$, where $\CC ^{\rm fd}$ is the class of $L_b$-structures of finite diameter,
and obtain the conclusions of Theorem \ref{T543} for the structures of the form $\bcd _{i<n}\X_i$, where $\X _i\in\CC ^{\rm fd}$.
\begin{te}\label{T520}
Let $\X =\la X,\r\ra =\bcd _{i<n}\X _i$, where $\X _i$, $i<n$, are pairwise disjoint $L_b$-structures of finite diameter, say $\leq m$,
let $\CT:=\Th (\X )$ and $\CT _i:=\Th (\X _i)$, for $i<n$. Then we have
\begin{itemize}\itemsep -0.4mm
\item[\rm (a)] $X /\r _{\rm rst}=\{ X_i:i<n\}$ and  $\r _{\rm rst}=D_{\ve _m ,\X}$ (see (\ref{EQ581}));
\item[\rm (b)] $\Y \equiv \bcd _{i<n}\X _i$ iff $\,\Y =\bcd _{i<n}\Y _i$, where $\Y _i \equiv \X _i$, for $i<n$;
\item[\rm (c)] $\BE \preccurlyeq \bcd _{i<n}\X _i$ iff $\;\BE =\bcd _{i<n}\BE _i$, where $\,\BE _i \preccurlyeq \X _i$, for $i<n$.
\item[\rm (d)] $\max \{I(\CT_i ):i<n \}\leq I(\CT )\leq \prod _{i<n}I(\CT_i ) $; thus, $I(\CT )=1$  iff $\,I(\CT_i)=1$, for each $i<n$;
\item[\rm (e)] If  the theories $\CT_i$, $i<n$, satisfy Vaught's conjecture, then $\CT$ satisfies Vaught's conjecture:
               $I(\CT )=\k:=\prod _{i<n}I(\CT _i)$, if $\k\in \{1,\o,\c \}$, and $I(\CT ) \in [3,\o )$, otherwise;
\item[\rm (f)] $\CT$ is atomic iff $\;\CT _i$, $i<n$, are atomic;
               then $\bcd _{i<n}\X_i^{\rm at}$ is a countable atomic model of $\CT$, where $\X_i^{\rm at}$ is a countable atomic model of $\CT _i$, for $i<n$;
\item[\rm (g)] $\CT$ is small iff $\;\CT _i$, $i<n$, are small;
               then $\bcd _{i<n}\X_i^{\rm sat}$ is a countably saturated model of $\CT$, where $\X _i ^{\rm sat}$ is a countably saturated model of $\CT_i$, for $i<n$.
\end{itemize}
\end{te}
\dok
(a) For $i<n$ we have $\delta (\X _i)\leq m$, so $\X _i$ is a connected substructure of $\X$.
Since $\la z,z'\ra \not\in \r$, whenever $z\in X_i$, $z'\in X_j$ and $i\neq j$,
the sets $X_i$, $i<n$, are maximal connected parts of $\X$,
which gives $X /\r _{\rm rst}=\{ X_i:i<n\}$.
If $X\models \ve _m [x,y]$,
then there is a path of length $\leq m$ from $x$ to $y$
and, hence, $\la x,y\ra\in \r _{\rm rst}$.
Conversely, if $\la x,y\ra\in \r _{\rm rst}$,
then $x,y\in X_i$, for some $i<n$,
and, since $\delta (\X _i)\leq m$,
there is a path of length $\leq m$ from $x$ to $y$ in $\X_i$ and, clearly, in $\X$,
which gives $\X\models \ve _m [x,y]$.
So, $D_{\ve _m, \X}=\r _{\rm rst}$.

(b) If $\Y\equiv \X$,
then by (a) and Proposition \ref{T519}(c) there is an enumeration $Y/D _{\ve _m ,\Y}=\{Y_i:i<n\}$ such that $\Y_i \equiv \X _i$, for each $i<n$.
Since $\X =\bcd _{i<n}\X _i$, for the sentence $\t :=\forall u,v \;(\neg \ve _m (u,v) \Rightarrow \neg R(u,v)) $ we have $\X \models \t$ and, hence, $\Y \models \t$.
If $i<j<n$, $y_i \in Y_i$ and $y_j \in Y_j$,
then, since $\Y \models \t$ and $\Y \models \neg \ve _m [y_i,y_j]$, we have $\la y_i ,y_j\ra \not\in R^\Y$;
thus,  $\Y = \bcd _{i<n}\Y _i$.
The converse is Fact \ref{T535}(b).

(c)  If $\BE \preccurlyeq \bcd _{i<n}\X _i$,
then by (a) $X /D_{\ve _m ,\X}=\{ X_i:i<n\}$
and, by Proposition \ref{T519}(b), $E /D_{\ve _m ,\BE }=\{ E\cap X_i:i<n\}$ and $\BE _i :=\BE \cap \X _i \preccurlyeq \X _i$, for $i<n$;
clearly we have $\BE =\bcd _{i<n}\BE _i$.
The implication ``$\Leftarrow$" is Fact \ref{T535}(c).

(d)--(g) These claims will follow from Theorem \ref{T543}, when we show that $\Perf _s (\bcd ,\CC ^{\rm fd})$.
First, since $|\bcd _{i<n}X_i|=\sum _{i<n}|X_i|$, we have (\ref{EQ621}).
Second, if a structure $\X\in\CC ^{\rm fd}$ is of diameter $\leq n$ and $\Y \equiv \X$,
then $\Y \models \f _{\d \leq n}$ and, hence, $\Y\in\CC ^{\rm fd}$;
thus the class  $\CC ^{\rm fd}$ is closed under $\equiv$.
Third, taking pairwise disjoint $\X _i, \Y _i\in\CC ^{\rm fd}$, $i\in \o$, we confirm (\ref{EQ615})--(\ref{EQ618}).
By (c) and Fact \ref{T535}(d) we have (\ref{EQ617}) and (\ref{EQ618}).
The implication ``$\Leftarrow$" in (\ref{EQ615}) follows from Fact \ref{T535}(a).
If $\X :=\bcd _{i<n}\X_i\cong \bcd _{i<m}\Y_i =:\Y$,
then $\Y \equiv \X$ and, by (b), $\Y$ has $n$ connected components;
thus $m=n$.
By Fact \ref{T535}(a) there is $\pi\in \Sym (n)$ such that $\X_i \cong \Y _{\pi (i)}$, for each $i<n$,
and the implication ``$\Rightarrow$" in (\ref{EQ615}) is true.
The implication ``$\Rightarrow$" in (\ref{EQ616}) follows from (b).
Conversely, if $\Y \cong \bcd _{i<n}\Y_i $, where $\Y _i \equiv \X _i$, for $i<n$,
then by Fact \ref{T535}(b) we have $\bcd _{i<n}\Y_i\equiv \bcd _{i<n}\X_i$
and, hence, $\Y\equiv \bcd _{i<n}\X_i$.
So, the implication ``$\Leftarrow$" in (\ref{EQ616}) is true.
\kdok
\section{The closure $\la\CC ^{\rm rt} \ra _{\du\Pi}$}\label{S6}
Here we extend the results of the previous two sections
to the minimal closure $\la \CC^{\rm rt} \ra _{\du\Pi}$ of the class $\CC^{\rm rt}$ under finite products and finite disjoint unions.
Generally speaking, if $\CC $ is an $\cong$-closed class of $L_b$-structures and $\X _i\in \CC$, for $i<n$,
then $\bcd _{i<n}\X _i$ is determined up to isomorphism, since there are $\X _i'\cong \X _i$, for $i<n$, belonging to $\CC$ and having pairwise disjoint domains. 
Also, the closures $\la \CC\ra _{\Pi}$ and $\la \CC\ra_{\du}$ of $\CC$ are obtained in one step
(because a finite product (union) of finite products (unions) is a finite product (union)),
while $\la \CC \ra _{\du\Pi} =\bigcup _{n\in \o}\CC _n$, where the classes $\CC _n$, $n\in \o$, are defined by recursion in a natural way.
But the following description (see Lemma 3.1(c) of \cite{Ksharp}) is more convenient.
\begin{fac}\label{T218}
If $\CC \subset \Mod _{L_b}$ is an $\cong$-closed class, then $\X\in \la \CC \ra _{\du\Pi}$ iff $\X$ is a finite disjoint union of finite products of structures from $\CC$.
\end{fac}
\begin{te}\label{T547}
For each partial order $\X\in\la \CC ^{\rm rt}\ra_{\du\Pi}$ we have
\begin{itemize}\itemsep -0.4mm
\item[\rm (a)] $\X = \bcd _{i<n}\prod _{j<m_i}\X _i^j$, where $n,m_i \in \N$ and $\X _i^j$, for $i<n$ and $j<m_i$, are rooted trees;
\item[\rm (b)] $\Y \equiv \X $ iff $\;\Y \cong \bcd _{i<n}\prod _{j<m_i}\Y _i^j$, where $\Y _i^j\equiv \X _i^j$, for $i<n$ and $j<m_i$;
\item[\rm (c)] $\BE \preccurlyeq \X $ iff $\;\BE = \bcd _{i<n}\prod _{j<m_i}\BE _i^j$, where $\BE _i^j\preccurlyeq \X _i^j$, for $i<n$ and $j<m_i$;
\item[\rm (d)] Vaught's conjecture is true for $\CT:=\Th (\X)$: if $\CT _i ^j:=\Th (\X _i^j)$, for $i<n$ and $j<m_i$, and $\k:=\prod _{i<n}\prod _{j<m_i}I(\CT _i^j)$, then
               $I(\CT)=\k$, if $\k \in \{ 1,\o ,\c \}$, and $I(\CT)\in \o \setminus 3$, otherwise;
\item[\rm (e)] $\CT$ is atomic iff $\;\CT _i^j$, for $i<n$ and $j<m_i$, are atomic;
               then $\bcd _{i<n}\prod _{j<m_i}(\X _i^j)^{\rm at}$ is a countable atomic model of $\CT$,
               where $(\X _i^j)^{\rm at}$ is a countable atomic model of $\CT _i^j$, for $i<n$ and $j<m_i$;
\item[\rm (f)] $\CT$ is small iff $\;\CT _i^j$, for $i<n$ and $j<m_i$, are small;
               then $\bcd _{i<n}\prod _{j<m_i}(\X _i^j)^{\rm sat}$ is a countably saturated model of $\CT$,
               where $(\X _i^j)^{\rm sat}$ is a countably saturated model of $\CT_i^j$, for $i<n$ and $j<m_i$.
\end{itemize}
\end{te}
\dok
(a) This is Fact \ref{T218}.
So the posets $\X _i :=\prod _{j<m_i}\X _i^j$, $i<n$, are pairwise disjoint
and of diameter $\leq 2$ (they have a smallest element).

(b) Let $\Y \equiv \X$.
We have $\X = \bcd _{i<n}\X _i $
and, by Theorem \ref{T520}(b), $\Y =\bcd _{i<n}\Y _i$, where $\Y_i \equiv \X _i$, for $i<n$.
By Theorem \ref{T544}(d) for each $i<n$ we have $\Y _i \cong \prod _{j<m_i}\Y _i^j$, where $\Y_i^j \equiv \X _i^j$, for $j<m_i$,
and by Fact \ref{T535}(a) $\Y =\bcd _{i<n}\Y _i \cong \bcd _{i<n}\prod _{j<m_i}\Y _i^j$.
Conversely, if $\;\Y \cong \bcd _{i<n}\prod _{j<m_i}\Y _i^j$, where $\Y _i^j\equiv \X _i^j$, for $i<n$ and $j<m_i$,
then $\Y \cong \bcd _{i<n}\Y _i $, where $\Y_i:=\prod _{j<m_i}\Y _i^j$, for $i<n$,
by Fact \ref{T043}(c)  for $i<n$ we have  $\Y_i\equiv \prod _{j<m_i}\X _i^j$
and by Fact \ref{T535}(b) $\Y \cong \bcd _{i<n}\Y _i \equiv \bcd _{i<n}\prod _{j<m_i}\X _i^j=\X$,
which gives $\Y \equiv \X$.

(c) Let $\BE \preccurlyeq \X =\bcd _{i<n}\prod _{j<m_i}\X _i^j$.
By Theorem \ref{T520}(c) $\BE =\bcd _{i<n} \BE _i$, where for $i<n$ we have $\BE _i \preccurlyeq \X _i:=\prod _{j<m_i}\X _i^j$
and, by Theorem \ref{T544}(c),  $\BE _i = \prod _{j<m_i}\BE _i^j$, where $\BE _i^j\preccurlyeq \X _i^j$, for $i<n$ and $j<m_i$.
Thus,  $\;\BE = \bcd _{i<n}\prod _{j<m_i}\BE _i^j$.
The converse follows from Facts \ref{T043}(d) and  \ref{T535}(c).

(d) For $i<n$ let $\X _i :=\prod _{j<m_i}\X _i^j$ and $\CT _i :=\Th (\X _i)$.
If $I(\CT _i^j)=\c$, for some $i<n$ and $j<m_i$,
then by Theorem \ref{T544}(c) $I(\CT_i)=\c$
and by Theorem \ref{T520}(e) $I(\CT)=\c$.
Otherwise we have $I(\CT _i^j)\leq \o$, for all $i<n$ and $j<m_i$;
and, by Theorem \ref{T544}(c), $I(\CT_i)\leq \o$.
Now, if $I(\CT _i^j)=\o$, for some $i<n$ and $j<m_i$,
then by Theorem \ref{T544}(c) $I(\CT_i)=\o$
and by Theorem \ref{T520}(e) $I(\CT)=\o$.
Otherwise we have $I(\CT _i^j)<\o$, for all $i<n$ and $j<m_i$;
and, by Theorem \ref{T544}(c), $I(\CT_i)< \o$.
Now, if $I(\CT _i^j)=1$, for all $i<n$ and $j<m_i$,
then by Theorem \ref{T544}(c) $I(\CT_i)=1$, for all $i<n$,
and by Theorem \ref{T520}(e) $I(\CT)=1$.
The remaining case is when $I(\CT _i^j)>1$, for some $i<n$ and $j<m_i$;
then $1<I(\CT _i)<\o $, by Theorem \ref{T544}(c) we have $3\leq\prod _{i<n}I(\CT _i)<\o$
and by Theorem \ref{T520}(e) $I(\CT)\in [3,\o)$.

(e) The theory $\CT$ is atomic
iff (by Theorem \ref{T520}(f)) for each $i<n$ the theory $\CT _i$ is atomic
iff (by Theorem \ref{T544}(d)) for each $i<n$ and each $j<m_i$ the theory  $\;\CT _i^j$ is atomic.
If $(\X _i^j)^{\rm at}\models \CT _i^j$, for $i<n$ and $j<m_i$, are countable atomic models,
then by Theorem \ref{T544}(d) $\prod _{j<m_i}(\X _i^j)^{\rm at}\models \CT _i$, for $i<n$, are countable atomic models
and by Theorem \ref{T520}(f) $\bcd _{i<n}\prod _{j<m_i}(\X _i^j)^{\rm at}$ is a countable atomic model of $\CT$.
The proof of (f) is similar to the proof of (e).
\kdok
\begin{rem}\label{R002}\rm
Concerning Theorem \ref{T547}(b) we note that, strictly speaking and working in ZFC, the class $\la \CC ^{\rm rt}\ra_{\du\Pi}$ is not $\cong$-closed.
Moreover, although the class $\CC ^{\rm rt}$ is first-order definable (and, hence $\equiv$-closed),
its closure $\la \CC ^{\rm rt}\ra_{\Pi}$ is not closed under $\cong$.
For example, the square $\X :=\la \o\times \o, \leq ^\X\ra$ of the ordinal $\o$ is a product of two rooted trees
and if $f:\o\times \o\rightarrow \o$ is a bijection,
and $\Y=\la \o \leq ^\Y\ra$, where $m\leq ^\Y n$ iff $f^{-1}(m)\leq ^\X f^{-1}(n)$, for $m,n\in \o$,
then, clearly, $\Y\cong \X$;
but $\Y$ is neither a tree nor a direct product of two or more rooted trees,
because its domain $\o$ is not a set of $n$-tuples, for $n\geq 2$.
Namely, $\emptyset \in \o$ but an ordered pair $\la x,y\ra :=\{\{x\},\{x,y\}\}$ is a non-empty set
and the same holds for $n$-tuples. So, $\Y\not\in \la \CC ^{\rm rt}\ra_{\du\Pi}$.
In fact we have
$$
\la\la \CC ^{\rm rt}\ra_{\du\Pi}\ra _{\cong}=\la\la \CC ^{\rm rt}\ra_{\du\Pi}\ra _{\equiv};
$$
namely, the inclusion ``$\subset$" is evident
and, conversely, if $\Y\in\la\la \CC ^{\rm rt}\ra_{\du\Pi}\ra _{\equiv}$,
then by Fact \ref{T218} and Theorem \ref{T547}(b) $\Y \cong \X$, for some $\X \in \la \CC ^{\rm rt}\ra_{\du\Pi}$
and, hence, $\Y\in\la\la \CC ^{\rm rt}\ra_{\du\Pi}\ra _{\cong}$.
Clearly, the results of Theorem \ref{T547} are valid for the partial orders from the closure $\la\la \CC ^{\rm rt}\ra_{\du\Pi}\ra _{\equiv}$.
\end{rem}

{\footnotesize

\end{document}